\documentclass[12pt]{article}

\usepackage{amsfonts}
\def\g{\mathfrak{g}}
\def\a{\mathfrak{a}}

\def\h{\mathfrak{h}}

\def\k{\mathfrak{k}}

\def\p{\mathfrak{p}}

\def\m{\mathfrak{m}}

\def\u{\mathfrak{u}}

\def\D{\mathbb{D}}
\def\F{\mathbb{F}}

\def\N{\mathbb{N}}

\def\Z{\mathbb{Z}}
\def\R{\mathbb{R}}
\def\C{\mathbb{C}}
\def\D{\mathbb{D}}

\def\fd {\hspace{0.35cm} \raise
-0.5mm\hbox{$\blacksquare$}\\}

\def\qed{{\null\hfill\ \raise3pt\hbox{\framebox[0.1in]{}}}\break\null}
\newtheorem{theo}{Th\'eor\`eme}
\newtheorem{prop}{Proposition}
\newtheorem{lem}{Lemme}
\newtheorem{cor}{Corollaire}
\newtheorem{rem}{Remarque}
\newtheorem{defi}{D\'efinition}

\def\ste{\par\smallskip\noindent}

\def\dem{ {\em D\'emonstration: \ste }}

\def\beq{\begin{equation}}
\def\eeq{\end{equation}}

\def\\D{\mathcal D}

\def\M{\mathcal M}
\newenvironment{res}
               {\begin{equation}\begin{minipage}{0.85\textwidth}}
               {\end{minipage}\end{equation}}
\def\ber{\begin{res}}
\def\eer{\end{res}}

\catcode`\Ž=\active\def Ž{\'e}
\catcode`\=\active\def {\`e}
\catcode`\ˆ=\active\def ˆ{\`a}
\catcode`\=\active\def {\`u}
\catcode`\=\active\def {\^{e}}
\catcode`\"=\active\def "{\^{i}}
\catcode`\™=\active\def ™{\^{o}}
\catcode`\ž=\active\def ž{\^{u}}
\catcode`\=\active\def {\c{c}}

\def\opto#1{{\smash{\mathop{\longrightarrow}\limits^#1}}} 

\title{Vecteurs distributions $H$-invariants de reprŽsentations 
induites, pour un espace symŽtrique rŽductif $p$-adique $G/H$.}

\author{ Philippe Blanc, Patrick Delorme\\ dŽdiŽ ˆ Alain Guichardet}
\date{}
\begin{document}
\maketitle
\noindent{\bf English title:}  $H$-fixed distribution vectors of induced representations,
for a reductive symmetric space $G/H$\ste \ste
\noindent{\bf RŽsumŽ: }Soit  $G$ le   groupe des points sur   $\F$ d'un groupe rŽductif
linŽaire dŽfini sur $\F$, un  corps local non archimŽdien de caractŽristique $0$. Soit
$\sigma$ une involution rationnelle de ce groupe algŽbrique dŽfinie sur
$\F$ et soit 
$H$ le groupe des points sur $\F$ d'un sous-groupe ouvert, dŽfini sur $\F$,  du groupe des
points fixes de
$\sigma$.  Nous construisons des familles de vecteurs  $H$-invariants 
dans le  dual de sŽries principales gŽnŽralisŽes, en utilisant l'homologie des groupes.
Des rŽsultats de   A.G.Helminck, S.P.Wang et 
A.G.Helminck, G.F.Helminck sur la  structure des espaces symŽtriques rŽductifs $p$-adiques 
sont aussi essentiels.  \ste \ste 
\noindent{\bf Summary:  }\ste 
Let $G$ be the group of $\F$-points of a linear reductive group defined over $\F$,  a non
archimedean local field of characteristic zero. Let 
$\sigma$ be  a rational involution of this group defined over $\F$ and let  $H$ be  the
group of $\F$-points of an open subgroup, defined over $\F$, of the group of  fixed points by 
$\sigma$. We built rational families of
$H$-fixed vectors in the dual of generalized principal series, using homology of groups.
Results of A.G.Helminck, S.P.Wang and A.G.Helminck, G.F.Helminck on the structure of
$p$-adic reductive symmetric spaces  are also essential.

\setcounter{section}{-1} \section{Introduction}

 Soit $G$ le groupe des points sur $\F$
d'un groupe linŽaire  algŽbrique rŽductif ${\underline G}$,  dŽfini sur $\F$, un 
  corps local non archimŽdien de caractŽristique $0$. Soit $\sigma$ une involution
de ce groupe algŽbrique 
${\underline G}$, dŽfinie sur
$\F$,
$H$ le groupe des points sur $\F$ d'un sous-groupe ouvert, dŽfini sur $\F$,  du groupe des
points fixes de
$\sigma$.
Cet article est destinŽ ˆ dŽbuter l'analyse harmonique sur l'espace symŽtrique rŽductif
$G/H$, en analogie avec celle sur les espaces symŽtriques rŽels (cf. [D] pour un 
survey sur ce sujet). 
\ste La premre Žtape que nous franchissons ici est la construction de familles de
formes linŽaires $H$-invariantes sur des sŽries principales gŽnŽralisŽes . Cette Žtape a ŽtŽ
franchie  dans le cas rŽel dans [BrD], ˆ l'aide de la thŽorie des
$D$-modules.
\ste  L'outil principal est ici l'homologie
lisse (voir plus bas). 
\ste Notre construction se limite ˆ celle de sous groupes paraboliques particuliers dits
$\sigma$-sous-groupes paraboliques. C'est l'expŽrience du cas rŽel qui conduit ˆ cette
restriction: dans ce cas ces familles de reprŽsentations suffisent ˆ dŽcrire la partie continue de la formule
de Plancherel pour l'espace symŽtrique. Un travail en cours d'Žlaboration de Nathalie
Lagier  prŽcisera le r™le jouŽ par les $\sigma$-sous-groupes paraboliques dans l'Žtude
des reprŽsentations admissibles de $G$ qui sont $H$-sphŽriques i.e. qui  possdent une
forme linŽaire non nulle $H$-invariante.  On notera, pour se fixer les idŽes,  que si 
$G$ n'admet pas de $\sigma$-sous-groupe parabolique diffŽrent de $G$ alors 
toutes les composantes isotropes de $G$ sont contenues dans $H$ ([HW], Lemme 4.5).\ste  
Il faut remarquer que   notre approche est  diffŽrente  de celle  de [Hi], [HiSat], [O]
 qui  dŽterminent, pour certains cas,   {\em }toutes les
reprŽsentations irrŽductibles ayant un vecteur non nul invariant par un bon  sous-groupe
compact maximal et une forme linŽaire non nulle
$H$-invariante, et obtiennent  une formule de Plancherel, ou bien   de l'approche  de [OS] qui
cherchent  toutes les reprŽsentations unitaires, irrŽductibles de $GL_{2n}$, ayant une forme
linŽaire non nulle invariante sous  le groupe symplectique. Dans [O] notamment  les fonctions
sphŽriques sont explicitŽes ˆ l'aide de polyn™mes de Macdonald, alors que dans  [OS], on
utilise des thŽormes sur les pŽriodes des formes automorphes.  
En particulier, nous ne cherchons pas ˆ expliciter compltement les zŽros et les p™les des
familles construites.  

\ste  On note $Rat  G$ le groupe des caractres rationnels de $G$ dŽfinis sur $\F$ et
$\a_G = Hom_{\Z}(RatG , \R )$  et on note $H_{G}:G\to \a_{G}$,
l'application qui,  ˆ $g\in G$,  associe l'application   $\chi \mapsto \vert
\chi(g)\vert_{\F}$. Alors
${\a}_G^* = Rat G \otimes_{\Z} \R$ et on dispose d'une application
surjective de
$({\a}_G^*)_{\C}$ dans l'ensemble $X(G)$ des caractres non ramifiŽs de $G$, qui ˆ
$\chi \otimes s$ associe le caractre $g \longmapsto \vert \chi(g) \vert_{\F}^s$. \ste 
L'involution $\sigma$ agit sur ${\a}_G$ et ${\a}_G^*$. On note ${\a}_{G,\sigma}$
(resp.
${\a}_{G}^{\sigma}$) l'ensemble des ŽlŽments de ${\a}_G$ antiinvariants
(resp. invariants) par
$\sigma$. Alors $\a_G = {\a}_{G,\sigma}\oplus {\a}_{G}^{\sigma}$ et ${{\a}_{G,\sigma}^*}$
s'identifie ˆ un sous-espace de $\a_G^*$. On note
$X(G)_{\sigma}$ l'image de $({{\a}_{G,\sigma}^*})_{\C}$ par l'application prŽcŽdente. Elle
possde une structure de variŽtŽ algŽbrique comme quotient d'un espace vectoriel par un
rŽseau.\ste On introduit des notations similaires  pour les sous-groupes de LŽvi
$\sigma$-stables de $G$.  On considre un sous-groupe parabolique $P$ de $G$ tel que
$\sigma(P)$ soit Žgal ˆ l'opposŽ de $P$ relativement ˆ un tore dŽployŽ maximal, $A_0$, 
contenu dans $P$ et $\sigma$-stable. On dit que $P$ est un $\sigma$-sous-groupe
parabolique. On note $M$ le sous-groupe de LŽvi $\sigma$-stable de $P$,  qui  est Žgal ˆ
$P \cap \sigma(P)$, et $U$ son radical unipotent. On note $A$ un tore dŽployŽ du centre
de $M$, maximal pour la propriŽtŽ d'tre contenu dans $\{ x \in G \vert
\> \sigma(x) = x^{-1}\}$.
\ste Alors ${\a}_{M,\sigma}^{*}$ s'identifie ˆ $\a^{*} = Rat A \otimes_{\Z} \R$. On note $X
= X(M)_{\sigma}$ et $B$ l'algbre des fonctions rŽgulires sur $X$.\ste 
Soit $(\delta, V_{\delta})$ une reprŽsentation lisse de $M$ irrŽductible, ou simplement
admissible de type fini, qu'on prolonge ˆ $P$ en la prenant triviale sur $U$. On introduit
pour
$\chi
\in X$, la reprŽsentation ${\delta}_{\chi} = \delta \otimes \chi$ de $M$. On dŽfinit
Žgalement une structure de $(M,B)$-module sur $V_{\delta} \otimes B$ en faisant agir $B$ par
multiplication sur le deuxime facteur et $m \in M$ par le produit tensoriel de $\delta(m)$
avec la multiplication par l'ŽlŽment 
$b_m$ de $B$ dŽfini par: $b_m(\chi) = \chi(m)$, $\chi \in X$.\ste
Soit encore:
$$\delta_B(m)(v\otimes b) = (\delta(m) v)\otimes b_mb$$ 
On Žtend l'action de $M$ ˆ $P$ en la prenant triviale sur $U$. On notera
$\delta_{\bullet}$ ˆ la place de  $\delta_{\chi}$ ou bien $\delta_B$.\ste 
 
On considre $ind_P^G V_{\delta_{\bullet}}$ l'ensemble des  $ \varphi: G \to
V_{\delta_{\bullet}}$ qui sont  invariantes ˆ gauche  par un
 sous-groupe  compact  ouvert  et telles que $$
\varphi(gmu) ={ \delta_{\bullet}} (m^{-1}) \varphi(g), \>\>  g \in G, m \in M, u \in U$$
On le note aussi  $I_{\bullet}$. Le groupe  $G$ agit sur $I_{\bullet}$ 
par la reprŽsentation rŽgulire gauche. De plus $B$ agit naturellement
sur $I_B$.\ste  L'objet de cet article est de dŽterminer les familles de
vecteurs $H$-invariants dans le dual de
$I_\chi$ dŽpendant, en un certain sens ( voir plus bas), de faon polynomiale de
$\chi$.\ste  Une remarque fondamentale pour notre travail est le fait que si $V$ est un
module lisse, l'espace des ŽlŽments $H$-invariants du dual de $V$, 
$V^{*H}$, est Žgal ˆ $H_0(H,V)^*$, o $H_0(H,V)$ est le quotient de $V$ par le
sous-espace engendrŽ par les vecteurs $gv - v$ o $g$ dŽcrit $H$ et $v$ dŽcrit $V$.\ste Dans
la premire partie de l'article, on Žtudie les foncteurs d'homologie lisse en degrŽ
supŽrieur. On les introduit gr\^{a}ce aux rŽsolutions projectives dans la catŽgorie des
modules lisses. Une sŽrie de propriŽtŽs sont brivement Žtablies(voir aussi [BlBr],
[Cas] et aussi [BoWal], [G] pour la cohomologie continue) notamment le
lemme de Shapiro, une utilisation des sous-groupes distinguŽs fermŽs qui
sont union  de sous-groupes compacts, la rŽsolution standard, le complexe
de chaines inhomognes. On Žtablit qu'une action naturelle du centre de
$G$ sur le complexe des chaines inhomognes induit une action triviale sur
l'homologie.\ste Tous ces rŽsultats sont
utilisŽs dans la preuve du ThŽorme principal. \ste  {\bf Hypothse simplificatrice,
pour l'introduction seulement}: On suppose que $HP$ est la seule $(H, P)$-double classe
ouverte dans
$G$.  \ste {\bf ThŽorme principal }:{\em
\ste (i) On note
$J_\chi$ = $\{
\varphi \in I_\chi \vert \> \varphi\hbox{  est ˆ support contenu dans } HP \}$ qui est un
sous-$H$-module lisse. Alors $H_0(H,J_\chi)$ est naturellement isomorphe ˆ $H_0(M \cap H
, V_\delta )$.\ste 
(ii) Il existe un polynome $q \in B$, non nul tel que pour tout $\chi \in X$
vŽrifiant $q(\chi) \not= 0$, l'injection naturelle de $J_\chi$ dans $I_\chi$ dŽtermine
un isomorphisme: $$H_0(H,J_\chi) \simeq H_0(H,I_\chi)$$
(iii) Passant aux duaux dans (i) et  (ii), si $\chi \in X$ et $q(\chi) \not= 0$, on dispose
d'un isomorphisme: $$V_\delta^{* M \cap H} \to I_\chi^{* H}$$ qu'on note: 
$$\eta \mapsto j(P,\delta,\chi,\eta)$$ 
On peut rŽaliser les $I_\chi$ dans un espace fixe $I$, par restriction des fonctions de
$G$ ˆ un sous-groupe compact maximal $K$. Alors le polyn\^ome $q$ peut tre choisi de
telle sorte que pour tout $\eta \in V_\delta^{* M \cap H}$ et $\varphi \in I$, $\chi
\mapsto q(\chi)<j(P,\delta,\chi,\eta),\varphi>$, $\chi \in X$, soit un ŽlŽment
de $B$, i.e une fonction rŽgulire sur $X$.\ste} 
Donnons une idŽe de la dŽmonstration de ce thŽorme.\ste 
On commence par Žtudier $I_\chi$ comme $H$-module. D'abord il existe un nombre fini de
$(H,P)$-doubles classes $HxP$ et on note $\Omega$ un ensemble, qui contient $e$,  de
reprŽsentants de celles-ci. On introduit des ouverts $O_0=HP \subset O_1
\subset ...  \subset O_n = G$ de sorte que
$O_{i+1}
\backslash O_i = H x_{i+1} P$ puis on note $I_i = \{\varphi \in I_\chi \vert supp
\varphi \subset O_i \}$  de sorte que
$I_0 = J_\chi$ et
$\{0\} \subset I_0 \subset I_1 ... \subset I_n$. On montre (voir aussi [BZ], ThŽorme 5.2)
que:
$$I_i / I_{i-1}
\simeq ind_{H \cap {x_i P x_i^{-1}}}^H \delta_{\chi \vert H \cap x_i P x_i^{-1}}^{x_i}$$ pour
$i=0,...,n$, o
$\delta_{\chi}^{x_i}$ est la reprŽsentation de $x_iPx_i^{-1}$ dans $V$  dŽfinie par:
$$ \delta_{\chi}^{x_i}(x_ipx_i^{-1}) = \delta_{\chi}(p), p \in P$$   
Ici les fonctions de l'espace de l'induite sont ˆ support compact
modulo le sous groupe induisant et sont invariantes ˆ gauche par un
sous-groupe compact ouvert de $H$.
\ste  En particulier:  $$J_\chi \simeq ind_{H \cap P}^H
(V_{\delta_\chi})$$
 Gr\^{a}ce au lemme de Shapiro, $H_0(H,J_\chi) = H_0(H \cap P ,
V_{\delta_\chi})$. Mais  $H \cap P = H \cap M$  car $P$ est un
$\sigma$-sous-groupe parabolique de $G$ et
$\chi \in X(M)_\sigma$ est trivial sur $M \cap H$. D'o le point (i) du ThŽorme.\ste 
Maintenant on choisit $x = x_i$ avec $i>0$. On pose: $P'=xPx^{-1}$,  
$\delta'=\delta^x$,   $\chi'=\chi^x$,   $M'=xMx^{-1}$ etc.\ste 
On veut calculer $H_*(H,I_i/I_{i-1})$. D'aprs le Lemme de Shapiro, cet espace est
isomorphe ˆ $H_*(H \cap P' , V_{\delta'_{\chi'}})$. \ste 
Soit
$V'=(M'\cap \sigma(U')) (U'\cap \sigma (P'))$. Alors 
$V'
\cap (M'
\cap
\sigma(M')) =
\{e\}$ et $V'$ est un sous-groupe unipotent distinguŽ dans $H \cap P'$. Comme $V' \cap H$
est un sous-groupe distinguŽ de  $P' \cap H = ((M' \cap \sigma(M') \cap H)(V' \cap
H)$, qui est de plus la  rŽunion de sous-groupes compacts, les propriŽtŽs de l'homologie
lisse montrent que:$$H_*(H
\cap P',V_{\delta'_{\chi'}}) = H_*(H \cap M',H_0(V' \cap H,V_{\delta'})_{\chi'})$$ Ici
l'indice
$\chi'$ indique, dans le second membre, la tensorisation par le caractre
$\chi'_{
\vert _{H
\cap M'}}$. On note que l'action de $V'$ sur $V_{\delta'_{\chi'}}$ ne
dŽpend pas de $\chi'$, car $\chi'$,  Žtendu ˆ
$P'$ en le prenant trivial sur $U'$, est trivial sur $V'$ puisque celui-ci est rŽunion de
sous-groupes compacts ouverts.\ste 
On a $P' \cap \sigma(P') = (M' \cap \sigma(M'))V'$. Un argument qui est une adaptation
au cas $p$-adique d'une idŽe de [CarD], montre alors que $H_0(V' \cap
H,V_{\delta'}) \simeq H_0(M' \cap \sigma(U'),V_{\delta'})$. Finalement on a:  $$
H_*(H,I_i / I_{i-1}) \simeq H_*(H \cap M',H_0(M' \cap \sigma(U'),{V_{\delta'})_{\chi'}})$$
Mais $M' \cap \sigma(U')$ est le radical unipotent du sous-groupe parabolique de $M'$, 
$M'
\cap
\sigma(P')$,  de sous-groupe de LŽvi $M' \cap \sigma(M')$.\ste 
Donc $H_0(M' \cap \sigma(U'),V_{\delta'})$ est un $M' \cap \sigma(M')$-module admissible de
type fini. Ceci implique  qu'il existe des caractres $\chi_1, ..., \chi_p$ du centre
$Z'$ de
$M'
\cap
\sigma(M')$ tels que pour tout $z \in Z$, $(z-\chi_1(z))...(z-\chi_p(z))$ agisse
trivialement sur celui-ci. Mais  on peut trouver un ŽlŽment $z_0$ de $Z' \cap H$ tel que
l'application
$b_{z_0}: X' \to \C$, donnŽe par $\chi' \mapsto \chi'(z_0)$, ne soit pas
constante. Alors, tenant compte du fait que l'action de $z_0$, qui est dans le centre
de $H \cap M' = H \cap M' \cap \sigma(M')$,  sur $H_*(M \cap H, H_0(M \cap
\sigma(U'),V_{\delta'})_\chi)$ est triviale, on voit facilement que si 
$q'(\chi) = (1-\chi_1(z_0)\chi'(z_0))...(1-\chi_p(z_0)\chi'(z_{0}))$ est non nul, l'
homologie ci-dessus est nulle, donc aussi $H_*(H,I_i/I_{i-1})$.\ste  Un argument
de longue suite exacte donne le point (ii) du ThŽorme, en utilisant pour $q$ le
produit des $q'$ lorsque $i$ varie dans $\{1,...,n\}$. \ste 
Pour Žtablir les propriŽtŽs de rationalitŽ de (iii), on reprend l'Žtude
ci-dessus en remplacant $\delta_\chi$ par $\delta_B$.
\ste\ste Les racines du plus grand tore dŽployŽ du centre de $M$,  $A_M$,  dans l'algbre de
Lie de
$P$ s'identifient ˆ des  ŽlŽmŽnts de $\a_M^{*}$, dont on note $\rho_{P}$ la demi-somme comptŽe
avec les mutiplicitŽs. \ste On note
$C(G,P,\delta^*,\chi)$ l'espace des fonctions
$\psi$ sur
$G$, ˆ valeur dans le dual $V_{\delta}^*$ de $V_{\delta}$, telles que:\steÊ\ber   Pour tout
$v \in V_{\delta}$, $g
\mapsto <\psi(g),v>$ est  continue  et :  $$ \psi(gmu) = e^{-2\rho_P(H_M(m))}
\chi(m) \delta^*(m^{-1}) \psi(g), g \in G, m \in M  , u \in U$$   \eer 
Le groupe $G$ agit par reprŽsentation rŽgulire gauche sur cet espace.\ste 
Si $\psi \in C(G, P,\delta^*,\chi)$ et $\varphi \in I_{\chi}$, on note
$<\psi,\varphi> = \int_K<\psi(k),\varphi(k)>dk$  \, qui dŽfinit un crochet de dualitŽ
$G$-invariant sur ces espaces. Ceci permet d'identifier $ C(G, P,\delta^*,\chi)$ ˆ un
sous-espace de $ I_{\delta_\chi}^*$.
\ste  On introduit une application $\varepsilon(G,P,\delta,\chi,\eta)$:  $G \to
V_{\delta^*}$,  caractŽrisŽe par la propriŽtŽ de covariance (0.1), nulle en dehors de $HP$,
$H$-invariante ˆ gauche et valant $\eta$ en $e$. \ste 
On introduit enfin une fonction   $\Vert\>\> \Vert$ sur $G/H$ et $M/{M \cap
H}$.\ste    \noindent {\bf ThŽorme }:\ste{\em Soit $\eta \in V_\delta^{*M \cap H}$
et
$r
\in
\R$ tel que pour tout
$v
\in V_\delta$, l'application de $M/M\cap H$ dans $\C$: $$m( M \cap H) \mapsto
<\delta^*(m)\eta,v>$$ soit bornŽe par un multiple de $\Vert m(M \cap H)
\Vert^r$. \ste  Alors il existe $\nu_0 \in \a_{M,\sigma}^* \subset \a_M^*$, $P$-dominant,
tel que pour tout $\chi \in X$ avec $Re\>\chi- 2\rho_{P}- \nu_0$ strictement
dominant par rapport aux racines de $A_M$ dans l'algbre de Lie de $P$,
la fonction $\varepsilon(G,P,\delta, \chi,\eta)$ est un ŽlŽment de
$C(G,P,\delta^*,\chi)$.\ste  L'ŽlŽment de $I_{\chi}^*$ correspondant ˆ
$\varepsilon(G,P,\delta, \chi,\eta)$  est Žgal ˆ
$j(P,\delta,\chi,\eta)$.}\ste 
On remarquera que si $P$ est un $\sigma$-sous-groupe parabolique minimal, tout $\eta \in
V_\delta^{*M \cap H}$ vŽrifie la condition du thŽorme. Des rŽsultats rŽcents de Nathalie
Lagier indiquent qu'il doit en tre ainsi mme si $P$ n'est pas minimal, au moins si $\delta$
est unitaire.
\ste  On notera que lorsque $G = G_1 \times G_1$ et $\sigma$ est l'involution qui Žchange les
facteurs, notre article donne une nouvelle approche des intŽgrales d'entrelacement (cf.[BrD]
pour le cas rŽel). 

Il est important de noter que les rŽsultats de A.G. Helminck et S.P. Wang [HWan], et  A.G.
Helminck et G.F. Helminck [HH], sur la structure des espaces
symŽtriques
$p$-adiques sont ˆ la base de ce travail. Il faut Žgalement dire que la page web de
A. G. Helminck annonce,  dans un rapport,  un travail avec G.F. Helminck sur le mme sujet
que celui traitŽ dans cet article. \ste  Nous espŽrons que nos rŽsultats seront
complŽmentaires  et que notre mŽthode, ˆ savoir l'utilisation de
l'homologie lisse,  pourra tre utile dans d'autres situations.\ste  {\bf
Remerciements}: \ste  Nous remercions Jean-Pierre Labesse pour avoir
rŽpondu ˆ de nombreuses questions pendant l'Žlaboration de ce travail. 

\setcounter{section}{0}

\section{Homologie lisse de $l$-groupes}
\setcounter{equation}{0}
\subsection{Modules projectifs}

On appelle $l$-groupe, un groupe localement compact, $G$, dŽnombrable ˆ l'infini, qui admet
une base de voisinages de l'ŽlŽment neutre formŽe de sous-groupes ouverts compacts. On
notera $dg$ une mesure de Haar ˆ gauche sur $G$.  On note
$C_c^\infty(G)$ l'espace des fonctions localement constantes ˆ support compact,
ˆ valeurs complexes. On note $\M(G)$ la catŽgorie des $G$-modules lisses, i.e.
des reprŽsentations de $G$ dans un espace vectoriel complexe dont tout vecteur
est fixŽ par un sous-groupe compact ouvert.
On voit facilement que  $C_c^\infty(G)$ muni de l'action rŽgulire droite, (resp.
gauche), est un $G$-module lisse. Appliquant ce rŽsultat ˆ $G^n$, puis plongeant
$G$ diagonalement dans $G^n$, on voit que $C_c^\infty(G^n)$  muni de l'action
rŽgulire gauche $L$, (resp. droite $R$), de $G$ (o $G$ agit diagonalement
sur $G^n$) est un $G$-module lisse.
Nous redonnons ici une preuve d'un rŽsultat de [Cas], ThŽorme A.4, qui est la version
$p$-adique d'un rŽsultat de [Bl] pour les groupes de Lie. 
\begin{lem}Le module  $C_c^\infty(G^{n+1})$ muni de l'action rŽgulire droite 
est projectif dans
$\M(G)$.
\end{lem}                                     
\dem On fixe ${\varphi} \in
C_c^\infty(G^{n+1})$ tel que 
$\int_{G^{n+1}}\varphi(x)dx = 1$.  
Si $g, \>x \in G^{n+1}, \>f \in C_c^\infty(G^{n+1})$,  on  dŽfinit:\beq f_g(x) =
f(x) {\varphi}(xg) \>\> \eeq
Alors $f_{g}$ est ŽlŽment de $C_c^\infty(G^{n+1})$.
Par un calcul immŽdiat on voit que:
\beq R_h f_{h^{-1}g} = (R_{h}f)_g  \eeq
Soit $P:  U \rightarrow V \rightarrow 0  $   un morphisme surjectif de
$G$-modules lisses et
     $F:  C_c^\infty(G^{n+1}) \rightarrow V $ un morphisme de
$G$-modules. Soit
$S$ une section linŽaire de la surjection $ P: U \rightarrow V$ \ste
On dŽfinit:$$
     {\overline F}(f):= \int_{G^{n+1}}{g_0Sg_0}^{-1} F(f_g) dg,\>\>  f
\in  C_c^\infty(G^{n+1})$$
Ici $ g = (g_0,...,g_n)$.
\ste On va vŽrifier  que ${\overline F}$ est un morphisme de $G$-modules entre
$C_c^\infty(G^{n+1})$ et $U$ tel que  $P\circ{\overline F} = F$.\ste
  Soit  
 $ f \in C_c^\infty(G^{n+1})$. Alors:   $$P({\overline F(f))} =
\int_{G^{n+1}}P(g_0Sg_0^{-1} F(f_g))dg$$        Mais $Pg_0 = g_0P$ et
$PS = Id_V$. Donc:  $$P\circ {\overline F}(f) =
\int_{G^{n+1}}F(f_g)dg$$   Mais 
$$\int_{G^{n+1}}F(f_g)dg = F(\int_{G^{n+1}}f_g dg)$$    car l'intŽgration est en
fait une somme finie. De plus   $$\int_{G^{n+1}}f_g(x)dg =
\int_{G^{n+1}}f(x)\varphi(xg)dg =f(x)$$   Il en rŽsulte que  $P{\overline F}(f) =
F(f) $.  Il faut maintenant voir que ${\overline F}$ est un $G$-morphisme.
Gr\^ace ˆ (1.2), on voit que pour $h = (h_0,...h_0) \in{G^{n+1}}$ avec $ h_0 
\in G$:   $${\overline F}(R_hf) = \int_{G^{n+1}}g_0Sg_0^{-1} F(R_hf_{h^{-1}g})
dg$$On effectue le changement de variable    $g' = h^{-1}g$,   qui implique que $
g_0 = h_0g'_0$ et l'on trouve      $${\overline F}(R_hf) =
\int_{G^{n+1}}h_0g_0Sg_0^{-1}F(f_g)dg = h_0{\overline F}(f)$$  \qed

 Il rŽsulte de ce qui prŽcde que si $V$ est un $G$-module lisse,
$C_c^\infty(G)\otimes V$ est projectif pour l'action rŽgulire droite sur le
premier facteur.
Remarquons que l'on peut remplacer l'action rŽgulire droite par l'action
rŽgulire gauche sur $C_c^\infty(G)$  car l'application   $f\rightarrow
f^\vee$  o     $f^\vee(g) = f(g^{-1})$    entrelace ces deux actions.
Montrons que: \ber {\em L'application de
$C_c^\infty(G)\otimes V$ dans $V$ dŽfinie par: $$f\otimes v \rightarrow
\int_Gf(g)gvdg$$ est un morphisme surjectif de  $G$-modules lisses pour le produit
tensoriel de la reprŽsentation rŽgulire gauche sur $C_c^\infty(G)$ avec la
reprŽsentation triviale sur $V$}\eer Il est surjectif car, si
$v$ est fixŽ par
$K$, on prend 
$f = 1_K/vol(K)$, o $1_{K} $ est l'indicatrice de $K$  et
$f\otimes v$ a pour image $v$. Comme $L_hf\otimes v$ a pour image
 $\int_Gf(h^{-1}g)gvdg$, on fait le changement de variable $g' = h^{-1}g$
 pour voir que c'est un $G$-morphisme. Donc $\M(G)$ a suffisamment de
projectifs et: 
\ber {\em Tout $V\in \M(G)$ admet une rŽsolution projective  par des modules du
type 
$C_c^\infty(G)\otimes W$, o $G$ agit par produit tensoriel de la
reprŽsentation rŽgulire gauche avec la reprŽsentation triviale sur $W$.}\eer 
Mais
on dispose d'une autre action de $G$ sur $C_c^\infty(G)\otimes V$:
 le
produit tensoriel des actions sur $C_c^\infty(G)$ et sur $V$. On remarque
que
$C_c^\infty(G)\otimes V$ est Žgal  ˆ $C_c^\infty(G,  V)$.
 Plus gŽnŽralement $G$ agit sur $C_c^\infty(G^{n+1},V)$
par:  
  \beq \label{action.}       (g.f)(g_0,...,g_n) = 
gf(g^{-1}g_0,...,g^{-1}g_n), \>\>f  \in
C_c^\infty(G^{n+1}, V),  \>\>g,g_0,..., g_n\in G \eeq  
 Notons que $C_c^\infty(G^{n+1}, V)$ est
isomorphe au $G$-module $C_c^\infty(G^{n+1})\otimes V$ o $G$ agit par la
reprŽsentation rŽgulire gauche de manire diagonale sur le premier
facteur et trivialement sur le second, l'isomorphisme Žtant donnŽ
par: \ber   $ f \mapsto \tilde{f}$  o $ \tilde{f} (g_0,...,g_n) =
g_0^{-1} f(g_0,...,g_n) $, $f \in C_c^\infty(G^{n+1}, V)$\eer Notons que
l'inverse de cet isomorphisme est donnŽ par: \beq  f \rightarrow
\widehat{f},
f
\in
C_c^{\infty}(G^{n+1}, V) \>\> avec:\>\> \widehat{f}(g_0,...,g_n) = g_0f(g_0,...,g_n)\eeq
car on a: 
  \ber $\tilde {(gf)}(x) =x_0^{-1}(gf)(x) = (x_0^{-1}g)f(g^{-1}x) =
(gx_0)^{-1}f(g^{-1}x) = \tilde{f}(g^{-1}x)$  \eer Donc les modules
$C_c^\infty(G^{n+1}, V)$ sont projectifs.
\subsection { Homologie lisse, rŽsolution standard et complexe de chaines
inhomognes  }
\begin{defi} Pour un module lisse
$V$, on dŽfinit
$V_G=J_G(V)$,  l'espace des coinvariants, comme le quotient de $V$ par l'espace vectoriel
engendrŽ par les
$gv - v$, $g \in G, v\in V$. Alors $J_G$ est un foncteur exact ˆ droite.
On dŽfinit alors l'homologie lisse de $V$, $H_*(G,V)$, comme Žtant l'homologie  du complexe:$$\cdots  \rightarrow J_G(P_1) \rightarrow J_G(P_0)
\rightarrow 0
$$  qui est obtenu par application du foncteur $J_{G}$ ˆ une rŽsolution projective de $V$
$$\cdots
\rightarrow P_2
\rightarrow P_1
\rightarrow P_0
\rightarrow V
\rightarrow 0 $$   Par des arguments
classiques (cf e.g. [CartE], Ch V, \S 3 ), $H_*(G,V)$ ne dŽpend pas de
la rŽsolution projective choisie.  De plus, $H_0(G,V)$ est canoniquement
isomorphe ˆ $J_G(V)$. \end{defi} En particulier on dispose de la rŽsolution
standard de $V$ (cf.  e.g. [G], Ch. I, Equation (12.2), pour le
cas des groupes discrets):
\begin{lem}  On considre la suite \beq \label{standard}\cdots\to  C_c^\infty(G^{n+2},V)
\opto{{\delta_n}} C_c^{\infty}(G^{n+1},V)
\opto{{\delta_{n-1}}}\cdots  \opto{{\delta_{0}}}C_c^\infty(G,V) \opto{\eta}  V\rightarrow 0
\eeq o:\beq  (\delta_n f)(g_0,...,g_n) = \sum_{i = o}^{n+1} (-1)^i\int_G
f(g_0,...,g_{i-1}, g,g_{i},...,g_n)dg \eeq
Ici:   \beq  \eta(f) =
\int_Gf(g)dg,\>   si f \in C_c^\infty(G, V) \eeq et l'action sur
$C_c^\infty(G^{n+1},V)$ est celle dŽfinie en (1.5). 
C'est une rŽsolution projective de $V$, dite rŽsolution standard.\end{lem}
\dem  Pour voir que la
suite (\ref{standard}) est exacte, on montre que pour chaque sous-groupe ouvert compact, la
suite:\beq  \to  C_c^\infty((K\backslash G)^{n+2},
V^K)\opto{{\delta_n}}C_c^\infty((K\backslash G)^{n+1}, V^K) \rightarrow \cdots
C_c^\infty(K\backslash G, V^K)
\opto{\eta} V^K \rightarrow 0 \eeq  est exacte. Pour cela on introduit: 
\ber$$ \sigma_n:
 C_c^\infty((K\backslash G)^n,V^K)\rightarrow C_c^\infty((K\backslash G)^{n+1},
V^K)
$$ 
$$(\sigma_{n}f)(g_0,...,g_n)=(1_K(g_0)/vol(K)) f(g_1,...,g_n)$$  $$avec 
\>\>(g_0,...,g_n)
\in G^{n+1}
\> et \> f\in C_c^\infty((K\backslash G)^n,V^K)$$  \eer
Alors $\sigma_n$ est une homotopie contractante comme on le vŽrifie
aisŽment,  car on a:$$\eta \sigma_0 = Id_{V^K}$$ et  
$$\delta_n\sigma_{n+1} +
\sigma_{n}\delta_{n-1} = Id_{C_{c}^{\infty}((K\backslash G)^{n+1}, V^K)}$$
puisque, pour $f\in C_{c}^{\infty}((K\backslash G)^{n+1}, V^K)$:  
$$(\delta_n\sigma_{n+1}f)(g_0,...,g_n) = (\delta_n((1_K/vol K)\otimes
f))(g_0,...,g_n)$$ 
$$\>\>\>\>\>\>\>\>\>\>\>\>\>\>\>\>\>\>\>\>\>\>\>\>
\>\>\>\>\>\>\>\>\>\>\>\>\>\>\>\>\>\>\>\>\>\>\>\>\>\>\>\>\>\>=f(g_0,...g_n)
- (1_K(g_0)/ vol(K))\otimes (\delta_{n-1}f)(g_1,...,g_n)$$ soit encore:
$$(\delta_n\sigma_{n+1}f)(g_0,...,g_n)=f(g_0,...,g_n) -((\sigma_n
\delta_{n-1})(f))(g_0,...,g_n)$$
qui est l'ŽgalitŽ voulue.
  Comme  tout $f \in C_c^\infty(G^{n+1}, V)$ est ŽlŽment de $
C_c^\infty((K/G)^{n+1}, V^K)$ pour $K$ sous-groupe compact ouvert assez
petit, l'existence de l'homotopie contractante $\sigma_n$ prouve
l'exactitude du complexe (\ref{standard}). Ainsi, tenant compte du Lemme
1 et de (1.5),  (\ref{standard}) est une rŽsolution projective de $V$.
\qed Donc $H_n(G,V)$ est le
$n$-{ime}
 groupe d'homologie du complexe: \beq\label{complexeinh}\cdots \opto{{
\delta'_{1}}} H_{0}(G, C_c^\infty(G^2, V))
\opto{{ \delta'_0}}  H_{0}(G, C_c^\infty(G,V)) \rightarrow 0 \eeq 
o $\delta'_n$ est dŽduit de $\delta_n $ par passage aux coinvariants. \ste DŽfinissons:
\ber \label{defT} $$T_n: C_c^\infty(G^{n+1}, V) \to
C_c^\infty(G^n,V)$$ par $$(T_n(f))(g_1,...,g_n) =
\int_G g^{-1}f(g,gg_1,...,gg_1...g_n)dg$$\eer 
\begin{lem}  Alors $T_n$ est surjectif et $Ker\>T_n$ est engendrŽ par les
ŽlŽments de la forme $g.f - f$, o on utilise les notations de (1.5).
\end{lem} Avant de procŽder ˆ la dŽmonstration du Lemme 3,
Žtablissons un autre Lemme.   
\begin{lem} (i) Soit
$\varphi
\in
C_c^\infty(G,W)$  o $W$  est un espace vectoriel. Alors si: $$
\int_G\varphi(g)d g = 0 $$ $\varphi$  est une combinaison linŽaire d 'ŽlŽments de la forme
$L_g\psi - \psi$, $\psi \in C_c^\infty(G,W), g \in G$.\ste 
(ii) L'application de $C_c^\infty(G,W)$ dans $W$ donnŽe par l'intŽgration sur
$G$ par rapport ˆ $dg$ passe au quotient en un isomorphisme entre
$H_{0}(G,C_c^\infty(G,W))$ et $W$.\end{lem}\dem (i) Comme
$C_c^\infty(G, W) = C_c^\infty(G)\otimes W$, et que $\varphi \in
C_c^\infty(G, W) $ ne prend qu'un nombre fini de valeurs, on se ramne au cas
o $W$ est de dimension finie, puis, en prenant une base de $W$, au cas
o $W$ est de dimension 1, i.e. $C_c^\infty(G, W) = C_c^\infty(G)$, ce
que l'on suppose dans la suite.  Soit donc $\varphi$ comme dans
l'ŽnoncŽ.\ste Alors
$\varphi$ est invariante ˆ droite par un sous-groupe compact ouvert assez
petit. \ste Donc  $$\varphi = \sum_{i=1}^n {\lambda_i} 1_{x_{i}K}$$  o $1_{xK}$ est l'
indicatrice de
$xK$.\ste
La condition sur $\varphi$ montre que  $\sum_{i=1}^n\lambda_i = 0$. Donc: 
$$\varphi =
\sum_{i=1}^n \lambda_i (1_{{x_i}K} - 1_K) = \sum_{i=1}^n {\lambda}_i(L_{x_i}1_K -1_K) $$
ce qui achve de  prouver (i). \ste Prouvons (ii). L'application est injective
d'aprs (i) et surjective d'aprs (1.3).\qed 

\begin{cor} Soit  $W$  une reprŽsentation lisse de $G$. On fait   agir $G$
sur 
$C_c^\infty(G,W)$ par la
reprŽsentation $$ (g,\varphi) \rightarrow (g.\varphi)  $$ en posant:
$$(g.\varphi)(x):=
g\varphi(g^{-1}x),  \>g,x \in G $$
  Si:
$$\int_{G} g^{-1}\varphi(g)dg = 0$$ alors $\varphi$ est combinaison linŽaire de fonctions de
la forme 
$g.\psi - \psi$ o $\psi \in C_c^\infty(G,W)$.  \end{cor} {\em DŽmonstration du Lemme 3:}
\ste D'abord, un simple changement de variable montre que $Ker\> T_n$ contient les
combinaisons linŽaires de fonctions du type $g.f-f$. Soit
$f
\in Ker \>T_n$. On applique  le corollaire ˆ $f$ regardŽe comme ŽlŽment de
$C_c^\infty(G,C_c^\infty(G^n,V))$, o l'on fait agir $G$ sur $W:= C_c^\infty(G^n,V)$ par:
$$(g\diamond h)(g_1,...,g_n):= gh(g^{-1}g_1,...,g^{-1}g_n), \>\>  
h \in C_c^\infty(G^n,V)$$  Alors $f$ est combinaison linŽaire de fonctions de la forme:
$g*F - F$, o $F \in C_c^\infty(G,W)$, $g \in G$  et: $$(g*F)(x):=
(g\diamond F)(g^{-1}x), \>\> x \in G$$ Alors, identifiant $F$ ˆ un
ŽlŽment de $C_c^\infty(G^{n+1}, V)$: $$(g*F)(g_0,...,g_n) =
gF(g^{-1}g_0,g^{-1}g_1,...,g^{-1}g_n)$$i.e., avec les notations de (1.5):
$$g*F = g.F$$ Ce qui prŽcde achve la dŽtermination du noyau de
$T_n$.\ste En utilisant (1.3), on voit de mme que si on
pose:$$(\widetilde{T}_nf)(g_1,...g_n) =
\int_Gg^{-1}f(g,gg_1,...,gg_n)dg$$ l'image de $\widetilde{T}_n$ est Žgale ˆ $W =
C_c^\infty(G^n,V)$.
Montrons que cela implique  que l'image de $T_n$ est aussi  Žgale ˆ
$C_c^\infty(G^n,V)$. \ste  Si
$ F 
\in C_c^\infty(G^n,V) $, on dŽfinit:
$$H(x_1,...x_n) = F(x_1,x_1^{-1}x_2,...,x_{n-1}^{-1}x_n), \>\> x_1,...,x_n\in G$$  donc $H
\in C_c^\infty(G^{n},V)$. De plus   
\beq \label{F=H}F(g_1,...g_n) = H(g_1,g_1g_2,...,g_1...g_n), \>\> g_1,...,g_n\in G \eeq 
\ste D'aprs ce qui prŽcde, il existe $f
 \in C_c^\infty(G^{n+1},V)$ tel que $\widetilde{T}_nf=F$.  Alors, tenant
compte de (\ref{F=H}) et (1.15), on voit que ceci implique que 
 $T_nf =  F$ et
ceci montre que $T_{n}$ est surjective.\qed 
\ste D'aprs le Lemme 3, $H_{0}(G, C_c^\infty(G^{n+1},V))$ est isomorphe
par
$T_n$ ˆ $C_c^\infty(G^{n},V))$. Un calcul direct montre  que:
\beq d_n(T_{n+1}(f)) = T_{n}(\delta_n(f)), \> f\in C_{c}^{\infty}(G^{n+2}, V)\eeq
o
  pour tout $n 
 \geq{1}$ et $  f\in C_{c}^{\infty}(G^{n+1}, V)$   \ber  $$(d_nf)(g_1,...,g_{n}) =
\int_Gg^{-1}f(g,g_1,...,g_{n})dg +$$ $$ \sum_{i = 1}^{n}(-1)^i
\int_Gf(g_1,...,g_{i-1},g,g^{-1}g_i,g_{i+1},...,g_{n_1})dg + (-1)^{n+1}
\int_Gf(g_1,...,g_{n},g)dg$$ \ste et
$$d_0 f = \int_Gg^{-1}f(g)dg - \int_Gf(g)dg, f\in C_c^{\infty}(G, V)$$\eer  
Donc \begin{prop} La diffŽrentielle
$\delta'_{n}$ transportŽe par ces  isomorphisme est Žgale ˆ $d_n$ et 
l'homologie lisse de
$V$ s'identifie ˆ l'homologie du complexe, dit complexe de chaines inhomognes:
\beq\cdots \to  C_{c}^{\infty}(G^{n+1}, V)\opto{{d_{n}}}C_{c}^{\infty}(G^{n},
V)\opto{{d_{n-1}}}\cdots \opto{{d_{0}}}V\rightarrow 0 \eeq  
\end{prop}
%[OMEttre en Si o: $f = 1_Kx \otimes {v}$\  o\ 
%$v\in V^K$. $g = kx$ \ste  $$d_0 f=\int_Gg^{-1} 1_{xK}(g)\otimes {v}\  
%dg  - 
%\int_G 1_{Kx}(g)
%\otimes {v}\   dg $$ 
%En faisant le changement de variable $g'=gx$, on trouve:  $$d_0 f=
%\delta_G(x)\int_{G} 1_{xKx^{-1}}(k) dk\otimes (x^{-1}v-v)$$   Ici
%$\delta_G$ dŽsigne la fonction modulaire de $G$].
%\ste 
\subsection{ Action de $G$ sur $H_{*}(G,V)$}
Pour tout
$g\in G$   on dŽfinit un automorphisme $R_g$ de $C_c^\infty(G^{n+1},V)$ muni
de la structure de 
$G$-module dŽfini en (\ref{action.}). L'opŽrateur $R_g$  correspond ˆ
l'action rŽgulire droite de la diagonale  i.e.: 
$$(R_g f)(g_0,...,g_n) = f(g_0g,...,g_ng),\>\>f\in  C_c^\infty(G^{n+1},V) $$

On vŽrifie sans peine que,  si $\Delta_{G} (g)= 1$, o $\Delta_{G}$ est la fonction
modulaire de $G$, $R_{g}$ est un automorphisme de $G$-modules de la rŽsolution standard
 de
$V$,  en complŽtant les flches $R_{g}$ par l'identitŽ de $V$. Donc il est homotope ˆ
l'automorphisme identitŽ, par une homotopie formŽe de morphismes de $G$-modules
(cf. e.g., [G] Chapitre I, Proposition 2.2).  Cette homotopie dŽtermine une homotopie du
complexe d'espaces vectoriels  (1.14) dont la cohomologie est l'homologie
lisse. Donc le
 morphisme $R_{g}$
induit l'identitŽ sur l'homologie lisse de $V$. 

\begin {prop}  \label{actionG} Pour tout $g \in G$ tel que  $\Delta_{G} (g)= 1$ et
tout $f \in C_c^\infty(G^n,V)$ on note: $$(\widetilde{R}_gf)(g_0,...,g_n) =
gf(g^{-1}g_0g,g^{-1}g_1g,...,g^{-1}g_ng)$$  
(i) Cette formule dŽfinit un automorphisme $\widetilde{R}_g$ du
complexe des chaines inhomognes (1.19), homotope ˆ l'identitŽ.\ste (ii) En particulier si
$f$ est un ŽlŽment du noyau de $d_n$, 
$ Z_n(G,V)$, alors 
 $ \widetilde{R}_gf - f $ est ŽlŽment de l'image de $d_{n+1}$,  $ B_n(G,V)$. \end{prop}
\dem
 Un calcul immŽdiat montre que $\widetilde{R}_g$, est le transportŽ de $R_g$ par $T_{n}$,
au scalaire multiplicatif $\Delta_{G}(g)$ prs.
\begin{cor} Si $z \in Z(G)$ est tel que $z-Id_V$ est inversible, $H_n(G,V)
=\{ 0\}$ pour tout
$n \in \N$. \end{cor} \dem    En effet, si $f \in Z_n(G,V)$, on vient de
voir que  
$\widetilde{R}_zf-f \in B_n(G,V)$.  Mais  
$((\widetilde{R}_z-Id)f)(g_1,...g_n) = (z-Id_V)f(g_1,...g_n)$. Il est
 alors facile de voir que $f \in B_n(G,V) $,  puisque
$(z-Id_V)^{-1}$ prŽserve $ B_n(G,V)$. Donc $H_n(G,V) =\{ 0\}$.
\qed
\subsection{Longue suite exacte d'homologie}
Le lemme suivant est standard. 
\begin{prop} Soit $0\to V\opto{i}W \opto{p} V/W\to 0$, une suite exacte courte
de $G$-modules lisses. Alors notons $i_{*}: H_{*}(G, V)\to H_{*}(G, W)$
( resp. $p_{*}: H_{*}(G, W)\to H_{*}(G,V/ W)$) les morphismes en homologie
dŽduits de $i$ et $p$. Pour $n\in \N^{*}$, il existe des morphismes canoniques
$c_{n}: H_{n}(G, V/W)\to H_{n-1}(G, V)$ tels que l'on ait la suite exacte
longue:
$$\cdots H_{n+1}(G, V/W)\opto{{c_{n+1}}} H_{n}(G, V)\opto{{i_{n}}} H_{n}(G, W) 
\opto{{p_{n}}} H_{n}(G,V/ W)\opto{{c_{n}}}\cdots \opto{{p_{0}}} H_{0}(G,V/ W)\to
0$$
\end{prop}
\subsection {Homologie lisse et sous-groupes distinguŽs} 
\begin{lem} \label{distingue}(i)  Si le $l$-groupe $G$ est une rŽunion de sous-groupes
compacts, le foncteur $H_0(G,.)$ est un foncteur exact sur ${\cal M}(G)$.
\ste (ii) De plus $H_i(G,V) = 0$ pour tout $i>0$ et tout  $V$ objet de ${\cal M}(G)$.
\end{lem}
\dem (i) est classique et (ii) rŽsulte de (i) et de la dŽfinition de l'homologie lisse.
\qed  
\begin{prop}  Soit $G$ un $l$-groupe, $H$ un sous-groupe fermŽ et distinguŽ de $G$ qui
est rŽunion de sous-groupes compacts ( donc unimodulaire). Alors si $V$ est un $G$-module
lisse, $H_0(H,V)$ est un $G/H$-module lisse et $H_*(G,V)$ est naturellement isomorphe ˆ
$H_*(G/H,H_0(H,V))$. \end{prop}
\dem L'assertion sur $H_0(H,V)$ est claire. D'aprs (1.4), le $G$-module lisse $V$ admet
une rŽsolution projective: $\cdots\to P_n \to P_{n-1} \to\cdots  \to P_0 \to V \to 0$,
telle que pour tout $n$, $P_n$ soit de la forme $C_c^\infty(G) \otimes W_n$ o $G$ agit par
le produit tensoriel de  la reprŽsentaion rŽgulire gauche sur $C_c^\infty(G)$ avec la
reprŽsentation triviale sur $W_n$. Alors
$H_0(H,P_n)
\simeq H_0(H,C_c^\infty(G)) \otimes W_n$ comme $G/H$-module. Mais $H_0(H,C_c^\infty(H))
\simeq
\C$, d'aprs le Lemme 4. 
Tenant compte de l'isomorphisme du Lemme  \ref{fibre} pour les actions rŽgulires gauches, on
en dŽduit que $H_0(H,C_c^\infty(G))$ s'identifie ˆ $C_c^\infty(H \setminus G)$ ( =
$C_c^\infty(G/H)$ puisque $H$ est distinguŽ dans $G$). On vŽrifie facilement que c'est un
isomorphisme de $G/H$-modules. Il en rŽsulte que $H_0(H,P_n)$ est un $G/H$-module lisse
projectif. Des  Lemmes 1 et 5 appliquŽs ˆ $H$, on dŽduit que: $$\cdots\to H_0(H,P_n) \to
\cdots
\to H_0(H,P_0) \to H_0(H,V) \to 0$$ est une rŽsolution projective du $G/H$-module lisse
$H_0(H,V)$. \ste 
Donc $H_*(G/H,H_0(H,V))$ est l'homologie du complexe:$$\cdots\to H_0(G/H,H_0(H,P_n)) \to
\cdots
\to H_0(G/H,H_0(H,P_0)) \to 0        $$
Mais il est clair que pour tout $G$-module lisse $W$, $H_0(G/H,H_0(H,W))$ s'identifie ˆ
$H_0(G,W)$. Alors le complexe dont l'homologie donne
$H_*(G/H,H_0(H,W))$ s'identifie ˆ celui dont l'homologie donne $H_*(G,V)$. D'o la
Proposition.\qed 

\subsection {ReprŽsentations induites comme espaces d'homologie }
  Soit $H$ un sous-groupe fermŽ du $l$-groupe $G$, et $(V,\pi)$ une reprŽsentation lisse
de $H$. On note $ind_{H}^{G}\pi$ la reprŽsentation rŽgulire gauche de $G$ dans l'espace
$ind_{H}^{G}V$ des fonctions $\varphi$ de $G$ dans $V$, invariantes ˆ gauche par un
sous-groupe compact ouvert de $G$,  ˆ support compact modulo $H$, et telles
que:$$\varphi(gh)  =
\pi(h)^{-1}\varphi(g), \>\> g \in G,\> h \in H$$ C'est une reprŽsentation
lisse de $G$.
\begin{prop}:  Si $f \in C_c^\infty(G)$  et $v \in V$, on dŽfinit:
$$i({f}\otimes{v})(g)
 = \int_{H}f(gh)\pi(h)vdh,  \>\> g \in G$$ Ici $dh$ est une mesure de Haar
 ˆ gauche sur $H$. \ste
(i) Alors  $i(f\otimes{v}) \in ind_{H}^{G}V$ et $i$ se prolonge en un
morphisme linŽaire de
$G$-modules entre $ C_c^\infty(G)\otimes V$,  muni du produit tensoriel de la reprŽsentation
rŽgulire gauche de
$G$ avec la reprŽsentation triviale sur $V$, et 
$(ind_H^GV, ind_H^G\pi)$, notŽ
$i$ ou
$i_{G,\pi}$.\ste (ii) Cette application est surjective. \ste 
(iii) Son noyau est Žgal ˆ l'espace engendrŽ par les $((R_h \otimes \pi(h)
\Delta_H^{-1}(h))\varphi) - \varphi$, $h \in H$, $\varphi\in
ind_{H}^{G}V$  o $\space R$ est la reprŽsentation rŽgulire droite de G
sur $C_c^\infty(G)$ et
$\Delta_H$ est la fonction modulaire de $H$. Celle-ci est caractŽrisŽe
par \beq
\int_{H}f(hh_{0})dh =
\Delta_H(h_{0})\int_{H}f(h)dh \hbox{ pour } h_0
\in H.\eeq (iv) Par passage au quotient $i$ dŽfinit un isomorphisme naturel
entre les
$G$-modules$H_0(H, C_c^\infty(G) \otimes V)$ et $(ind_H^{G}V, ind_H^{G}\pi)$,  o
$C_c^\infty(G) \otimes V$  est muni du produit tensoriel de l'action
rŽgulire droite de
$H$ avec $\Delta_{H}  \pi$. Ici $G$ agit sur $H_{0}(H,C_c^\infty(G) \otimes
V)$ par passage au quotient du produit tensoriel $L \otimes Id$ de la
reprŽsentation rŽgulire gauche sur $C_c^\infty(G)$ avec la reprŽsentation
triviale sur $V$.\ste   (v) Si $U$ est un ouvert $H$-invariant ˆ droite, on
notera
$(ind_H^{G}V)(U)$ ou
$C_c^\infty(U, \pi)$ l'espace formŽ des $\varphi \in ind_H^{G}V$ tels que le support de
$\varphi$ est contenu dans $U$. Alors $C_c^\infty(U, \pi)$ est l'image de $C_c^\infty(U)
\otimes V$ par $i$. 
\end{prop}
\dem  (i) rŽsulte d'un simple changement de variable. \ste
(ii) On va montrer (v) qui implique (ii).
Soit $\varphi \in C_c^\infty(U,\pi)$. Alors $\varphi$ est  invariante ˆ gauche par un
sous-groupe ouvert compact $K$, donc de support de la forme 
 $\bigcup_{x \in X}KxH$.
Les $K x H$ sont ouverts car $K$ est ouvert. \ste Le support de $\varphi$
Žtant compact modulo $H$, $\varphi$ est nulle en dehors de la rŽunion 
$\bigcup_{j}Kx_{j}H$,   pour un nombre fini de $x_{j}$, o 
l'union est disjointe et contenue dans $U$. On note $v_j =
\varphi(x_j)$. Les propriŽtŽs d'invariance ˆ gauche de $\varphi$ par $K$ et de
covariance ˆ droite sous $H$ montrent que
$v_{j}$ est invariant par $x_j^{-1}Kx_j \cap H$. \ste 
  L'intŽgrale $\int_{H}1_{{x_j}^{-1}Kx_j \cap H}(h)dh$ est non nulle car
$x_j^{-1}Kx_j \cap H$ est ouvert dans $H$. On note $c_j$  son inverse. Nous
affirmons que: $$\varphi = i(
\Sigma_jc_j(f_j \otimes v_j)), \hbox{ o  } f_j = 1_{Kx_j}$$
En effet les deux membres ont les mmes propriŽtŽs d'invariance ˆ  gauche sous
$K$ et de covariance ˆ droite sous $H$. Ils  sont nuls en dehors de
$\bigcup_j Kx_jH$, et sont
 dŽterminŽs
par leur valeurs en les $x_j$. Or la valeur en $x_j$ du second membre est
Žgale ˆ
$$c_j\int_H 1_{Kx_j}(x_jh) \pi(h) v_j dh = c_j\int_H 1_{{x_j}^{-1}Kx_j
\cap H}(h)
\pi(h)v_jdh = v_j$$ 
comme dŽsirŽ. Ceci achve de prouver (v) et (ii).\ste Prouvons  (iii).  On va
se ramener ˆ dŽterminer le noyau de $i_{G,\pi}$ dans le cas o $G = H$. En
effet, d'aprs l'appendice, Lemme \ref{fibre}, pour l'action rŽgulire droite  de $H$,
$C_c^\infty(G)$ est isomorphe ˆ $C_c^\infty(G/H)
\otimes C_c^\infty(H)$ o $H$ agit par la reprŽsentation rŽgulire droite 
sur le deuxime
facteur et trivialement sur le premier. Notons $s$ une section continue de la projection
$G\to G/H$ (cf. [M]). De mme, on a un isomorphisme $T$ d'espaces vectoriels entre
$C_c^\infty(G/H)
\otimes ind_H^{H}V$ et
$ind_H^{G}V$, dŽfini par:
$$
T(\psi \otimes f)(s(x)h) = \psi(x)f(h), \>\>\psi \in C_c^\infty(G/H), \>\>f
\in ind_H^{H}V,\>\> x\in G/H,\>\> h\in H$$ d'inverse: $$(T^{-1}(\varphi)(x,h) =
\varphi(s(x)h),\>\> x\in G/H,\>\> h\in H,\>\>
\varphi \in ind_H^{G}V$$Dans ces
isomorphismes $i_{G,
\pi }$ s'identifie ˆ 
$id_{C_c^\infty(G/H)} \otimes i_{H, \pi }$. La dŽtermination du noyau de
$i=i_{G,\pi}$, se rŽduit donc bien ˆ celle de $i_{H, \pi }$. On suppose
dŽsormais que $G = H$ et l'on note $\Delta$ au lieu de $\Delta_{H}$. 
\ste Soit $f\in i_{G, \pi}$. Alors $$\int_G \pi(g) f(g) dg = 0$$
Mais $\Delta(g) dg$ est une mesure de Haar invariante ˆ droite sur $G$, notŽe $d_r g$.\ste
Donc $$\int_G \pi(g) \Delta(g)^{-1} f(g) d_r g = 0$$
On applique ˆ la fonction sous le signe somme le Lemme 4 (ou plžtot sa version ˆ
"droite"). Alors $$\pi(g) \Delta(g)^{-1}f(g) = \Sigma_{i=1}^n (R_{x_i} \varphi_i)(g) -
\varphi_i(g), \>\> g\in G $$o  $\varphi_i \in C_c^\infty(G,V)$ et $x_i \in G$.\ste On pose,
pour $g\in G$, $\psi_i(g) = \pi(g)^{-1} \Delta(g) \varphi_i(g)$, de sorte que
$\varphi_i(g) = \pi(g) \Delta(g)^{-1} \psi_i(g)$. \ste
Alors$$f(g) = \Sigma_{i=1}^n (\pi(g)^{-1} \Delta(g) \pi(gx_i) \Delta(gx_i)^{-1}
\psi_i(gx_i) - \psi_i(g) , g\in G $$
D'o $$ f = \Sigma_{i=1}^n \pi(x_i) \Delta(x_i)^{-1} R_{x_i} \psi_i - \psi_i $$
comme dŽsirŽ.\qed

\subsection{Lemme de Shapiro en homologie lisse}
\begin{lem}
\label{shapiro} Si $(\pi,V)$ est un $G$-module lisse, $H_0(G,ind_H^G V)$ est
naturellement isomorphe ˆ
$H_0(H,V)$. Si $T: V \rightarrow V'$ est un morphisme de $H$-modules,
l'application entre $H_0(G,ind_H^G V)$ et $H_0(G,ind_H^G V')$,   dŽduite de l'application
entre  $H_0(H,V)$ et
$H_0(H,V')$  par  les isomorphismes ci-dessus,  est Žgale ˆ  l'application dŽduite du
morphisme induit, $ind  T$,  entre 
$ind_H^GV$ et $ind_H^GV'$, par passage au quotient.\end{lem}
\dem En effet il rŽsulte de la Proposition 5 (iv) que:$$ind_H^G V \simeq
H_0(H,C_c^{\infty}(G)
\otimes V)$$ o $C_c^\infty(G) \otimes V$ est muni du produit tensoriel de
l'action rŽgulire droite de $H$ avec $\Delta_H^{-1} \pi$. Cette action
commute avec l'action rŽgulire gauche de $G$. Donc: $$H_0(G,ind_H^GV) \simeq 
H_0(G \times H,C_c^\infty(G) \otimes V)$$ qui est aussi isomorphe ˆ
$H_0(H,H_0(G,C_c^\infty(G)\otimes V))$. Mais tenant compte du Lemme 4, on a
finalement:
$$H_0(G,ind_H^GV) \simeq  H_0(H,V)$$ o $H$ agit sur $V$ par $\pi$. 
L'assertion sur les morphismes est immŽdiate en suivant les isomorphismes
prŽcŽdents. \qed
\begin{prop} 
(Lemme de Shapiro pour l'homologie lisse)\ste  Soit $V$ un $H$-module lisse. Alors
$H_*(G,ind_H^GV)$ est naturellement isomorphe ˆ $H_*(H,V)$.\end{prop}
 {\em DŽbut de la dŽmonstration de la Proposition 6}\ste  Soit:  
$$\opto{{u_n}}P_n\cdots\opto{{u_0}}P_0 \to V \to0$$ 
une rŽsolution projective du $H$-module lisse $V$ par des modules projectifs
du type
$C_c^\infty(H) \otimes W_n$, o $H$ agit par la reprŽsentation rŽgulire
gauche sur
$C_c^\infty(H)$. \ste
 Alors $ind_H^G(P_n)$ est un $G$-module lisse projectif isomorphe ˆ $C_c^\infty(G)
\otimes W_n$  d'aprs le Lemme suivant. 
\begin{lem}
Soit $W$ un espace vectoriel. Le $G$-module $ind_H^GC_c^\infty(H)$ induit du
$H$-module $C_c^\infty(H)$ pour l'action rŽgulire gauche de $H$ est isomorphe ˆ
${C_c}^\infty(G)$ muni de l'action rŽgulire gauche de $G$.\end{lem}
\dem 
Clairement: $ind_{\{e\}}^H \C = C_c^\infty(H)$, o $\C$ est le $H$-module
trivial et
$C_c^\infty(H)$ est muni de l'action rŽgulire gauche. Utilisant l'induction par
Žtages, qui s'Žtablit facilement,  on en dŽduit le lemme.\qed
{\em Fin de la dŽmonstration de la Proposition 6}\ste
Donc
$$\cdots \opto{{ind\>u_n}}ind_H^G
P_n\cdots \to ind_H^G P_{n-1}\cdots\opto{{ind\>u_0}} ind_H^G P_{0} \to
ind_H^G V \to 0$$ est une rŽsolution projective du $G$-module lisse
$ind_H^G V$.\ste
L'homologie lisse $H_{*}(G,ind_H^G V)$ est naturellement isomorphe, d'aprs
le Lemme 6, ˆ l'homologie du complexe: $$\cdots\to H_0(H,P_n) \to \cdots \to H_0(H,P_0)
\to 0$$D'o l'isomorphisme voulu.\qed 
\subsection {Filtration d'une reprŽsentation induite} 
 On suppose maintenant que $P$ et $H$ sont deux sous-groupes fermŽs d'un  $l$-groupe
$G$ dŽnombrable ˆ l'infini et que $G$ n'admet qu'un nombre fini de $(H,P)$-doubles
classes. Utilisant l'action de $H \times P$ sur $G$, $H$ agissant ˆ
gauche et $P$ ˆ droite, on dŽduit de l'Appendice , Lemme  \ref{orbites}, qu'il existe des
ouverts: $$U_0 = \emptyset
\subset U_1 \cdots \subset U_n = G$$ 
tels que pour $i \geq 0$, $U_i\setminus U_{i-1}$ soit une double classe
$Hx_iP$  ouverte dans $G \setminus U_{i-1}$. 
\begin{prop} \label{filtration}
Soit $(\delta, V)$ une reprŽsentation lisse de $P$. On note:  $$I_i =
\{\varphi \in ind_P^G V\vert  \>\>le \>\>support\>\> de\>\> \varphi \>\>
est\>\> contenu \>\> dans \>\>
 U_i\}$$
\ste  Alors $I_i$ est $H$-invariant et la reprŽsentation de $H$ dans $I_{i}
 /
I_{i-1}$ est isomorphe ˆ
 $ind_{H \cap x_iPx_i^{-1}}^H \delta^{x_i}_ {\mid H \cap x_iPx_i^{-1}}$,
o
$\delta^{x_i}$ est la reprŽsentation de $x_iPx_i^{-1}$ dans $V$  dŽfinie par:
$$ \delta^{x_i}(x_ipx_i^{-1}) = \delta(p), p \in P$$ 
\end{prop}
\dem
En utilisant l'isomorphisme naturel de $ind_P^G\delta $ et $ind_{x_iP{x_i}^{-1}}^G
\delta^{x_i}$, on se ramne ˆ dŽmontrer l'assertion sur $I_i /
I_{i-1}$  dans le cas o $x_i$ est l'ŽlŽment neutre de $G$, ce que l'on
fait dans la suite. On notera $I$ au lieu de $I_i$, $J$ au lieu de
$I_{i-1}$, $U$ au lieu de $U_i$, $U'$ au lieu de $U_{i-1}$.\ste On associe ˆ  $f \in I$ sa
restriction ˆ $H$, $r_H(f)$.  Montrons que $\varphi:=r_H(f)$ est un ŽlŽment de l'espace
$ind_{H \cap P}^H$ 
$\delta_{\mid H \cap P}$.\ste La seule chose qui n'est pas immŽdiate est que
le support de $\varphi$ est compact modulo $H \cap P$. Le support de $
\varphi$ est Žgal ˆ
$Supp f
\cap H$.\ste  Soit $(h_n)$ une suite dans $Supp \varphi$. Il faut trouver
$(x_n)$ suite dans $H
\cap P$ telle que $(h_nx_n)$ ait une sous-suite convergente. Comme $Supp f$ est compact
modulo $P$, il existe une suite $(p_n)$ dans $P$ telle que $(h_np_n)$
admette une sous-suite convergente.\ste
 Notons $ H \times_{H \cap P} P$  le quotient de $H \times P $ par
 $H \cap P$, agissant ˆ droite sur le premier facteur et ˆ gauche sur le
second.
D'aprs l'Appendice,  Lemme \ref{orbites}, appliquŽ ˆ l'action de  $H \times P$ sur $G$,
l'application produit
$H
\times P \to HP$ passe au quotient en un isomorphisme topologique de $H
\times _{H \cap P} P$ sur $HP$.\ste  Il en rŽsulte qu'il existe une
suite $(x_n)$ dand $H \cap P$ telle que $(h_nx_n)$ et $(x_n^{-1}p_n)$
admettent des sous-suites convergentes. Ceci achve de prouver la
compacitŽ du support de $r_H(f)$ modulo $H \cap P$ et donc que $r_H(f)$ est
ŽlŽment de l'espace  $ind_{H \cap P}^H \delta_{\mid H \cap P}$.\ste  Le
reste de la dŽmonstration consiste ˆ montrer que $r_H$ est surjective, de
noyau
$J$.\ste 
Soit $f \in Ker \> r_H$. Cela signifie que $f$ est nulle sur $H$ donc sur $HP$.  
Elle est donc ˆ support dans $U \setminus HP$ = $U'$, i.e. $f \in J$. Il
reste seulement ˆ prouver la surjectivitŽ.\ste 
On note $r$ la restriction des ŽlŽments de $C_c^\infty(U) \otimes V$ ˆ
$C_c^\infty(HP) \otimes V$ qui est surjective car $HP = U \setminus U'$ est fermŽ
dans $U$.\ste
On dispose d'une application $i_{P,\delta}$ de $C_c^\infty(U_i) \otimes V$ dans $I_i$ qui
est surjective, (cf. Proposition 5). \ste 
Enfin on dŽfinit une application $j_{H,\delta}$ de $C_c^\infty(HP) \otimes V$
dans l'espace de $ind_{H \cap P}^H \delta_{\mid H \cap P}$, en posant:
$$(j_{H,\delta}(f))(h) = \int_P \delta(p)f(hp)dp, f \in C_c^\infty(HP)
\otimes V, h \in H$$ \ste 
o $dp$ est une mesure de Haar ˆ gauche sur$P$ (la mme que celle utilisŽe dans
la dŽfinition de $i_{P,\delta}$).\ste 
Alors on a: 
\beq   j_{H,\delta} \circ r = r_H \circ i_{H,\delta }  \eeq
Tenant compte de l'isomorphisme topologique entre $H \times_{H \cap P} P$ et
$HP$, on voit qu'on a une flche surjective $i_{H,P}$ de 
$C_c^\infty(H) \otimes C_c^\infty(P) \otimes V$ dans $C_c^\infty(HP) \otimes
V$, donnŽe par:$$i_{H,P}(\varphi \otimes \psi \otimes v)(hp) = \int_{H
\cap P}
\varphi(hx) \psi(x^{-1}p)dx$$ o $dx$ est une mesure de Haar ˆ gauche   sur
$H \cap P$.\ste Pour prouver la surjectivitŽ de $r_H$, tenant compte de celle de $r$ et de
(1.21), il suffit de dŽmontrer la surjectivitŽ de la composŽe: $$k:= j_{H,\delta} \circ
i_{H,P}$$  La dŽfinition de
$k$ comporte deux intŽgrations successives. D'aprs le  ThŽorme de Fubini, on
peut intervertir l'ordre des intŽgrations (les fonctions considŽrŽes sont ˆ
support compact). On trouve ainsi, pour $\varphi \in C_c^\infty(H), \psi \in
C_c^\infty(P)$ et $v \in V$ 
$$(k(\varphi \otimes \psi \otimes v))(h) =
\int_{H \cap P}(\int_P \varphi(hx)\psi(x^{-1}p) \delta(p) v dp) dx$$ On change alors $p$
en $x^{-1}p$ et on tient compte de l'invariance ˆ gauche  de la mesure $dp$: 
$$(k(\varphi \otimes \psi \otimes v )(h) = \int_{H \cap P} \varphi(hx) \delta(x) (\int_P
\psi(p) \delta(p) v dp) dx$$  On prend pour $\psi$ l'indicatrice d'un
sous-groupe compact ouvert de $P$, laissant fixe $v$, divisŽe par le volume de
ce sous-groupe.\ste Alors:
$$k(\varphi \otimes \psi \otimes v ) = i_{H,\delta_{\mid H \cap P}}(\varphi
\otimes v)$$  Tenant compte de la Proposition 5, ceci achve de montrer la
surjectivitŽ de $k$. Ceci achve la preuve de la Proposition.\qed 
 \section{Vecteurs distributions $H$-invariants de reprŽsentations
induites}
\setcounter{equation}{0}
\subsection{Groupes rŽductifs, notations}
On va utiliser largement des notations et conventions de [Wald].
Soit $\F$ un corps local non archimŽdien,
 de caractŽristique 0, de corps rŽsiduel $\F_q$. 
On considre divers groupes algŽbriques dŽfinis sur $\F$, et
on utilisera des abus de terminologie du type suivant: "soit $A$ un tore
dŽployŽ " signifiera " soit
$A$ le groupe des points sur $\F$ d'un tore dŽfini et dŽployŽ
sur
$\F$ ". Soit $G$ un groupe linŽaire algŽbrique rŽductif et
connexe.
 On fixe un sous-tore $A_{0}$ de $G$, dŽployŽ et maximal
pour cette propriŽtŽ, dont le centralisateur est notŽ $M_{0}$. Si $P$ est
un sous-groupe parabolique contenant
$A_{0}$, il possde un unique sous-groupe de LŽvi contenant $A_{0}$, notŽ
$M$ (notŽ aussi
$M_{P}$). Son radical unipotent  sera notŽ
$U$ ou
$U_{P}$. 
  Si $H$ est un groupe algŽbrique, on note $Rat (H)$ 
le groupe des caractres rationnels de $H$ dŽfinis sur $\F$. Si $E$
est un espace vectoriel , on note $E^*$ son dual. S'il est de plus rŽel, on note   $E_{\C}
$ son complexifiŽ. On note $A_{G}$ le plus grand tore dŽployŽ dans le centre
de $G$. \ste On note
$\a_{G}= Hom_{\Z} (Rat (G),\R)$. La restriction des caractres
rationnels de $G$ ˆ $A_{G}$  induit un isomorphisme: 
\beq  Rat(G)\otimes _{\Z}\R \simeq Rat(A_{G})\otimes _{\Z} \R\eeq  
On dispose de l'application canonique, $H_{G}: G \rightarrow \a_{G}$
 dŽfinie par: 
\beq \label{H} e^{<H_{G}(x), \chi>}= \vert \chi (x)\vert_{\F}, \> x\in G,
\chi
\in Rat (G)\eeq o
$\vert. \vert_{\F}$ est la valuation normalisŽe de $\F$. Le noyau de 
$H_{G}$, qui  est notŽ   $G^{1}$,  est l'intersection des noyaux des
caractres de $G$ de la forme $\vert \chi \vert_{\F}$, $\chi \in Rat
(G)$. On notera
$X(G)= Hom (G/G^{1}, {\C}^*)$.  On a des notations similaires pour
les sous-groupes de LŽvi. 
\ste Si $P$ est un sous-groupe
parabolique contenant $A$,  on notera
$\a_{P}= \a_{M_{P}}$, $H_{P}=H_{M_{P}}$. On note $\a_{0}= \a_{M_{0}}$,
$H_{0}=H_{M_{0}}$.
   On note $\a_{G, \F}$, resp. ${\tilde
\a}_{G, \F}$ l'image de $G$, resp. $A_{G}$, par $H_{G}$. Alors
$G/G^{1}$ est un rŽseau isomorphe ˆ
$\a_{G,\F}$.
Soit $M$ un
sous-groupe de LŽvi contenant $A_{0}$.  Les inclusions:
$A_{G}\subset A_{M}\subset M\subset G$, dŽterminent   un morphisme 
surjectif
 $\a_{M, \F}\rightarrow  \a_{G, \F}$, resp.  un morphisme injectif
 ${\tilde
\a}_{G, \F}   \rightarrow {\tilde
\a}_{M, \F}$, qui se prolonge de manire unique  en une  application
linŽaire surjective entre $\a_{M}$ et $\a_{G}$, resp. injective entre
$\a_{G}$ et
$\a_{M}$.
La deuxime  application permet d'identifier $\a_{G}$ ˆ un sous-espace de
$\a_{M}$ et le noyau de la premire, $\a^{G}_{M}$, vŽrifie ;
\beq \label{oplus} \a_{M}= \a^{G}_{M}\oplus \a_{G}\eeq 
 Il y a une surjection:
\beq  \label{surjection}(\a_{G}^{*})_{\C}\rightarrow X(G)\rightarrow
1\eeq qui est dŽfinie comme suit. En utilisant la dŽfinition   de $\a_{G}$, en remplaant $G$
par 
$A_{G}$, gr\^ace ˆ (2.1), on associe  ˆ
$\chi\otimes s$,  le caractre $g\mapsto \vert \chi (g)\vert
^{s}$. Le noyau est un rŽseau et ceci dŽfinit sur $X(G)$ une structure de
variŽtŽ algŽbrique complexe pour laquelle $X(G)\simeq \C^{*d}$, o
$d=dim_{\R}\a_{G}$. Pour $\chi \in X(G)$, soit $\nu \in
\a_{G,\C}^{*}$ un ŽlŽment se projetant sur $\chi $ par l'application
(2.4). La partie rŽelle $Re\>\nu \in \a_{G}^{*}$ est indŽpendante du
choix de $\nu$. Nous la noterons $Re \> \chi$. Si $\chi \in Hom (G,
\C^{*})$ est continu, le caractre $\vert \chi \vert$ appartient ˆ $X(G)$. On pose
$Re\> \chi=  Re\> \vert \chi\vert $. De mme, si $\chi \in Hom(A_{G},
\C^{*})$ est continu, le caractre $\vert \chi \vert$ se prolonge de faon unique en un
ŽlŽment de $X(G)$  ˆ valeurs dans $\R^{*+}$, que l'on note
encore  
$\vert
\chi \vert$  et on pose $Re \> \chi=  Re\> \vert \chi
\vert$.
\ste On dŽfinit 
$Im
\> X(G):=
\{
\chi
\in X(G)\vert Re
\>  \chi =0\}$ l'ensemble des ŽlŽments unitaires de $X(G)$. \ste De
l'isomorphisme naturel (2.1) on dŽduit aisŽment l'ŽgalitŽ:
\beq A_{G}^{1}= A_{G}\cap G^{1}\eeq On  voit facilement que $A_{G}^{1}
$ est le plus grand sous-groupe compact de $A_{G}$.
On note $X_{*}(G)$ l'ensemble des sous-groupes ˆ un paramtre de 
 $A_{G}$. C'est un rŽseau. On fixe une fois pour toute une
uniformisante. On note alors $\Lambda (G)$, l'image de $X_{*}(G)$ dans
$G$ par l'application '' Žvaluation  en l'uniformisante'', qui est un
rŽseau isomorphe ˆ $X_{*}(G)$ par cette Žvaluation. En effet, tout
ŽlŽment  de  $X_{*}(G)$  est dŽterminŽ par sa valeur sur une uniformisante
$\varpi$: en Žcrivant $A_{G}$ comme un produit de tores dŽployŽs de
dimension 1, on se ramne ˆ une assertion sur les sous groupes ˆ un
paramtre d'un tore dŽployŽ de dimension 1, qui est claire. Pour  tout
ŽlŽment non trivial de
$\Lambda(G)$, il existe un caractre non ramifiŽ rŽel de
$A_{G}$ non Žgal ˆ 1 sur celui-ci: on se ramne immŽdiatement aux tores
dŽployŽs de dimension 1. Ce caractre
non ramifiŽ de
$A_{G}$ ce prolonge en un caractre non ramifiŽ de $G$, d'aprs ce qui
prŽcde. Appliquant ceci aux sous-groupes de Levi, on a:
\ber  Pour tout ŽlŽment non trivial de $\Lambda (A_{M})$, il existe
$\chi \in X(M)$ diffŽrent de 1 sur cet ŽlŽment.\eer 
Soit  $A$ est un tore dŽployŽ de $G$, et $\lambda \in \Lambda(A)$. On note  $P_{\lambda}$ le
sous-groupe parabolique contenant $A$ pour lequel les racines
$\alpha$  de $A$ dans  l'algbre de Lie de $P_{\lambda}$ vŽrifient
$\vert \alpha(\lambda)\vert_{\F}$ est supŽrieur  ou Žgal ˆ 1.
Alors on a:\beq P_{\lambda}=\{ g\in G\vert 
(\lambda^{-n}g\lambda^{n})_{n\in \N} \hbox{  est    born\'ee}  \} = 
\{ g\in G\vert 
(\lambda^{-n}g\lambda^{n})_{n\in \N}  \hbox{ converge}\}\eeq
Alors:
\ber $M_{\lambda}:=\{g\in G\vert \lambda^{-1}g\lambda=g\}$ est le
sous-groupe de LŽvi de $P_{\lambda}$  contenant $A$. \eer
et \ber $U_{\lambda}:=\{ g\in G\vert 
(\lambda^{-n}g\lambda^{n})_{n\in \N}   \hbox{  converge   vers  
e}\}$ est le radical unipotent de $P_{\lambda}$.\eer
Si $A$ est un tore dŽployŽ maximal, tout sous-groupe parabolique de $G$
contenant
$A$ est de cette forme.\ste
Soit $\sigma$ une involution, dŽfinie sur $\F$,  du groupe algŽbrique dont $G$ est le groupe
des points sur
$\F$. Soit 
$H$ le groupe des points sur $\F$ d'un sous-groupe ouvert, dŽfini sur $\F$,  du groupe des
points fixes de
$\sigma$.\ste  On suppose maintenant que
$P$ est un sous-groupe parabolique quelconque de $G$. Alors (cf. [HWan], Lemme 2.4)  il
contient un tore dŽployŽ maximal de
$G$, $A$, qui est $\sigma$-stable.  Donc $P=P_{\lambda}$,  pour un ŽlŽment
$\lambda$ de $\Lambda(A)$.
Alors si $g \in P_{\lambda}$, $g$ = $mu$   avec $m \in
M_{\lambda}$, $u \in  U_{\lambda}$ et $m = lim_{n \to +\infty}
\lambda^{-n} g \lambda^n$. Enfin montrons que:  \ber \label{suite} Si 
$(g_n)$ est une suite dans
$P_{\lambda}$ qui converge vers $g$, $(\lambda^{-n}g_n \lambda^n)$ 
tend vers la limite de $(\lambda^{-n}g\lambda^{n})$ \eer
En effet,  on se rŽduit facilement ˆ $(g_n)$, suite
dans $U$. Mais, par rŽduction ˆ $GL(n)$, gr\^ace ˆ la linŽaritŽ de $G$,
 on voit que: \ber Si $X$ est une partie bornŽe de $U$  et $(x_n)$ une suite
dans $X$, alors $(\lambda^{-n}x_n\lambda^n)$ tend vers $e$ \eer Ceci
implique (\ref{suite}).
\subsection{Intersection d'un sous-groupe parabolique avec $H$} 
\ste
  On suppose maintenant que
$P$ est un sous-groupe parabolique. Alors
$P$ contient un tore dŽployŽ maximal de $G$, $A$, qui est $\sigma$-stable
(cf. [HWan], Lemme 2.4), et $P$ est de la forme $P_{\lambda}$ pour un
 $\lambda \in \Lambda (A)$. On note $M: = M_\lambda$, $U:= U_\lambda$,
 $\mu:= \lambda\sigma(\lambda) = \sigma(\lambda)\lambda$, $Q:=
P_{\mu}$, $V$ son radical unipotent, qu'on Žvitera de confondre avec l'espace d'une
reprŽsentation,  et
$L:= M_\mu$.\ste
\begin{prop}:\ste (i) On a:$$P\cap\sigma(P) = (M\cap \sigma(M))V'$$ o $$V' =
(P\cap\sigma(P)) \cap V$$  De plus $V'$ est distinguŽ dans
$P\cap\sigma(P)$ et $V'\cap ( M\cap\sigma(M)) = \{e\}$.\ste(ii) On a:$$V' =
(M\cap\sigma(U))(U\cap\sigma(P))$$ $$U\cap\sigma(P) = (U\cap\sigma(M))
(U\cap\sigma(U))$$ o $U\cap\sigma(P)$ est distinguŽ dans $V'$ et
$U\cap\sigma(U)$ est distinguŽ dans $U\cap\sigma(P)$.\ste(iii)
On a:$$P\cap H = (M\cap\sigma(M) \cap H)(V' \cap H)$$ o $
V' \cap H$ est distinguŽ dans $P\cap H$.\ste(iv) Pour tout $x \in
M\cap \sigma(U)$, il existe $h \in H\cap V'$ et $y  \in U\cap
\sigma(P)$ tels que $x = hy$.\end{prop}
\dem (i) Soit $g \in P\cap \sigma(P)$. Alors la suite
$(\lambda^{-n}g\lambda^n)$ est bornŽe. De plus comme $\sigma(A) =A$
on a $\lambda \in \sigma(P)$. Donc $(\lambda^{-n}g\lambda^n)$ est une
suite bornŽe dans $\sigma(P)$. Il en rŽsulte facilement que
$(\sigma(\lambda)^p \lambda^{-n} g \lambda^n
\sigma(\lambda)^p)_{n,p \in \N}$ est un ensemble bornŽ dans $\sigma(P)$. En
particulier, $(g_n) = (\sigma(\lambda)^{-n}\lambda^{-n} g \lambda^n
\sigma(\lambda)^n)$ est bornŽe et $g$ est ŽlŽment de $Q = LV$. Donc $g =
lv$ avec $l \in L$, $v \in V$, o $l$ est la limite de la suite
$(g_n)$. Comme $\lambda\in A$, $\sigma(\lambda) \in P\cap \sigma(P)$. 
 Mais $(\lambda^{-n} g \lambda^n)$ converge dans $P\cap \sigma(P)$,  en
particulier dans $\sigma(P)$,  vers $m
\in M$.\ste Alors d'aprs (\ref{suite}),
$(g_n)$ converge dans $\sigma(P)$ vers la mme limite que
$(\sigma(\lambda)^{-n} m  \sigma (\lambda)^n)$, qui est ŽlŽment de
$\sigma(M)$. Donc $l\in \sigma(M)$ et par raison de symŽtrie, $l
\in M \cap \sigma(M) \subset P\cap \sigma(M)$. Alors $v \in 
P \cap \sigma(P) \cap V = V'$. On a donc bien:$$P \cap \sigma(P) = (M
\cap \sigma(M))V'$$ Maintenant $M \cap \sigma(M) \subset L$ et $V'
\subset V$, donc leur intersection est rŽduite ˆ un ŽlŽment. Enfin $P
\cap \sigma(P)$ est contenu dans $Q$, $V'$ est contenu dans $V$ qui est
distinguŽ dans $Q$. Donc $V'$ est distinguŽ dans $P \cap \sigma(P)$.\ste  Soit
$x \in V' \subset P \cap \sigma(P)$. On Žcrit $x = mu$, avec
$m \in M$, $u \in U$. Ici $m$ est la limite de $(\lambda^{-n} x
\lambda^n)$ qui est une suite dans $P \cap \sigma(P)$ (car $\lambda \in
P\cap \sigma(P)$ puisque $\sigma(A)=A$), donc dans $P$. Il rŽsulte de (\ref{suite}) que la
limite de
$(\sigma(\lambda)^{-n}\lambda^{-n} x \lambda^n \sigma(\lambda)^n)$ est 
Žgale ˆ la limite de 
$(\sigma(\lambda)^{-n} m \sigma(\lambda)^n)$. Comme $x \in V$, cette
limite est $e$, donc $m \in \sigma(U)$. Finalement $m \in M \cap
\sigma(U)$. Par ailleurs $A$, donc  aussi  $\lambda$ et
$\sigma(\lambda)$ normalisent $Q$ et $V$, donc aussi $V'$. Ainsi, la
limite $m$ est ŽlŽment de $V'$, donc  $u$ aussi. En particulier $u \in
\sigma(P)$, et finalement $u \in U \cap \sigma(P)$.\ste Ceci prouve
l'inclusion $V' \subset (M \cap \sigma(U))(U \cap \sigma(P))$. L'inclusion
inverse est Žvidente. Ceci prouve la premire ŽgalitŽ de (ii).  \ste 
Maintenant soit
$u
\in U
\cap
\sigma(P)$. On Žcrit $x = m'u'$ avec $m'
\in \sigma (M)$, $u' \in \sigma(U)$. Alors $m'$ est la limite de
$(\sigma(\lambda)^{-n} x \sigma(\lambda)^n)$. Comme $x \in U$,
$(\lambda^{-n} x \lambda^n)$ tend vers $e$. Une application rŽpŽtŽe de
(\ref{suite}), montre que: $$lim_n(\lambda^{-n} m' \lambda^n) =
lim_n(\lambda^{-n}
 \sigma(\lambda)^{-n} x \sigma(\lambda)^n \lambda^n)
= lim_n(\sigma(\lambda)^{-n} {e}  \sigma(\lambda)^{n}) = e$$  Donc $m'
\in U \cap \sigma(M)$. \ste 
On montre de mme que $u' \in U\cap \sigma(U)$. Donc: $U\cap\sigma(P)
\subset (U \cap \sigma(M))(U \cap \sigma(U))$.   L'inclusion inverse
est Žvidente et ceci achve de prouver la deuxime ŽgalitŽ de (ii).\ste 
Comme $\sigma(U)$ est distinguŽ dans $\sigma(P)$, $P \cap \sigma(U)$ est
distinguŽ dans $V' \subset P \cap \sigma(P)$. On a donc prouvŽ (ii).\ste 
Montrons (iii).   Comme $\sigma( \mu ) = \mu$, $Q$ est $\sigma$-stable et
$\sigma(V) =V$. \ste  Soit alors $p \in P \cap H$. On Žcrit, gr\^ace ˆ (i),  
$p = mv$ avec $m \in M \cap \sigma(M)$, $v \in V' = V \cap P \cap
\sigma(P)$. Alors $p = \sigma(m) \sigma(v)$  avec $\sigma(m) \in 
M \cap \sigma(M)$, $\sigma(v) \in V'$.\ste  Mais ces deux
dŽcompositions doivent coincider: $\sigma(m) = m$ et $\sigma(v) = v$.\ste D'o
l'on dŽduit l'ŽgalitŽ de (iii).\ste  Le fait que $V' \cap H$ est distinguŽ
dans $P \cap H$ rŽsulte de (ii) et de l'ŽgalitŽ: \ste$P \cap H = \sigma(P) \cap
H$.  Ceci achve de prouver (iii).\ste 
(iv)  Notons $\m$ l'algbre de Lie de $M$,
$\sigma(\u)$ l'algbre de Lie de $\sigma(U)$. Le groupe
$L\cap\sigma(U)$ est un groupe algŽbrique unipotent, de m\^eme que $V'$. Son application
exponentielle est donc dŽfinie et  surjective.\ste 
Soit $x \in M \cap \sigma(U) \subset V'$. Alors $x = exp(X)$ avec $X \in
\m \cap \sigma(\u)$. Alors $X$+$\sigma(X)$ est ŽlŽment
de l'algbre de Lie de $V'$. On note $h = exp(X+\sigma(X))$, c'est un
ŽlŽment de $V' \cap H$. Alors $h$ et $x$ ont la mme projection dans
$V'/(U \cap \sigma(P)) \simeq M \cap \sigma(U)$ ˆ savoir: $exp(X)$. 
Donc $x = hg$ avec $g \in U \cap \sigma(P)$. Ceci achve la preuve de la
proposition. \qed \subsection{$H\cap P$-homologie lisse} Soit $X$ une
sous-variŽtŽ  de
$X(M)$, la variŽtŽ des caractres non ramifiŽs de $M$. On note $B_X$, l'algbre
des fonctions ˆ valeurs complexes sur $X$ engendrŽe par les fonctions
$b_m$, $m \in M$, o $b_m(\chi) = \chi(m)$, pour $\chi \in X$. On la 
notera souvent  $B$ au lieu de $B_{X}$. Alors si
$(\delta,V_{\delta})$ est un $M$-module lisse, $V_{\delta_{B}} = V_{\delta}
\otimes B$ est un $(M,B)$-module pour la reprŽsentation lisse $\delta_B$, de $M$,
dŽfinie par:\beq \label{deltaB}\delta_B(m)(v \otimes b) = \delta(m) v \otimes
{b_m b}, 
   \>\>b\in B, m \in M \eeq
et $B$ agissant par multiplication sur le deuxime facteur. On Žtend cette action de $M$ ˆ
$P$ en faisant agir $U$ trivialement.\ste 
\begin{lem} \label{HcapP} Avec les notations ci-dessus, $H_*(P \cap H,
V_{\delta_B})$ est un $B$-module naturellement isomorphe ˆ $H_*(M \cap H,
H_0(M \cap \sigma(U), V_{\delta})_B)$. Ici  $H_0(M \cap
\sigma(U),V_{\delta})$ est muni d'une structure  naturelle de $M \cap
\sigma(M)$-module, car $M\cap \sigma(U)$ est le radical unipotent du
sous-groupe parabolique de $M$, $M\cap \sigma(P)$, de sous-groupe de
LŽvi $M\cap \sigma(M)$.\end{lem} 
\begin{dem} Le groupe  $P \cap H$ admet pour sous-groupe distinguŽ $V' \cap H$,
avec  $M \cap H$ pour  quotient. De plus $V'$, donc aussi $V' \cap M$,
est rŽunion de sous-groupes compacts, comme sous-groupe fermŽ du
radical unipotent $V$ de $Q$. Donc (cf. Proposition 4),  $H_*(P \cap H,
V_{\delta_B})$  est isomorphe ˆ   $H_*(M \cap H, H_0(V' \cap H,
V_{\delta_B}))$. On vŽrifie aisŽment que c'est un isomorphisme de
$B$-modules.\ste  Il s'agit donc d'Žtudier $H_0(V' \cap H, V_{\delta_B})$. On
va dŽmontrer que pour toute reprŽsentation lisse $(\delta, V_{\delta})$ de
$M$, Žtendue ˆ $P$ en faisant agir $U$ trivialement, la surjection
naturelle de
$H_0(V' \cap H,  V_{\delta})$ dans $H_0(V', V_{\delta})$ est un
isomorphisme. \ste Il suffit pour cela de prouver la surjectivitŽ de
l'application transposŽe: $$H_0(V',V_{\delta})^*  \to H_0(V' \cap H,
V_{\delta})^*$$ Mais, il est clair que $H_0(V' \cap H, V_{\delta})^*$
est Žgal ˆ $Hom_{V' \cap H} (V_{\delta},\C)$\ste  Soit $T\in Hom_{V' \cap
H}( V_{\delta}, \C )$ et soit  $x \in M \cap \sigma(U)$. D'aprs la Proposition \ref{HcapP}
(v), on a:    
\ber $x = hy$ pour un  $h \in H \cap V'$ et un  $y \in U \cap \sigma(P)$\eer  
Si $v \in V_{\delta}$, on a:$$T(\delta(x)v) = T(\delta(hy)v)$$ et
finalement: $$T(\delta(x)v) = T(v)$$ car $\delta (y)v = v$ puisque
$y \in U $ et $T(\delta(h)v) = T(v)$ car $T \in Hom_{H \cap
V'}(V_{\delta}, \C)$.\ste 
Par ailleurs si $x \in U \cap \sigma(P)$, $\delta(x)v = v$ et, gr\^ace ˆ la
Proposition \ref{HcapP} (ii), on a  finalement $T \in Hom_{V'}(V_{\delta},\C)$. Ceci
prouve l'isomorphisme voulu.\ste
L'isomorphisme  $H_0(V' \cap H, \delta_B) \simeq H_0(V',\delta_B)$
est clairement un isomorphisme de $B$-modules. Mais comme les
caractres non ramifiŽs de $M$ sont triviaux sur le groupe unipotent $M
\cap \sigma(U)$,  $M \cap \sigma(U)$ agit sur $V_{\delta_B} =
V_{\delta} \otimes B$ par $\delta$ sur le premier facteur et
trivialement sur le second. Finalement on a un isomorphisme de
$B$-modules:$$H_0(M
\cap
\sigma(U), V_{\delta_B}) \simeq H_0(M \cap \sigma(U),
V_{\delta})\otimes B$$Par ailleurs $H_0(M \cap \sigma(U), V_{\delta})$ de
mme que $H_0(M \cap \sigma(U),  V_{\delta_B})$ est un $M \cap
\sigma(M)$-module, puisque $M \cap \sigma(M)$ est un sous-groupe de
LŽvi du sous-groupe parabolique $M \cap \sigma(P)$ de $M$. Alors
l'isomorphisme ci-dessus  dŽtermine  un isomorphisme de $(M\cap \sigma(M), B)$-modules
entre:
$$H_0(M
\cap
\sigma(U), V_{\delta_B}) \simeq H_0(M \cap \sigma(U),
V_{\delta})_B$$
Ceci achve la preuve du Lemme.\qed 

 \end{dem} 
\begin{prop}  \label{z} 
On suppose que $V_{\delta}$ est un $M$-module lisse admissible de type fini, que $X$ contient
le caractre trivial de $M$,  qu'il existe $z_0$ ŽlŽment de
l'intersection du centre de $M \cap \sigma(M)$ avec $H$ et qu'il existe  $\chi
\in  X$ tels que  $\chi(z_0) \neq 1$. Alors il existe $q \in B=  B_X$, non 
constant, produit d'ŽlŽments de la forme $1-cb_{z_0}, c\in
\C^*$,  tel que, pour tout $p\in
\N$,  le
$B$-module:
$H_{p}(H
\cap M, H_0(M
\cap
\sigma(U), V_{\delta})_{B})$ soit annulŽ par $q$.\end{prop} 
\dem Le $M \cap \sigma (M)$-module, $H_0(M \cap \sigma(U),
V_{\delta})$ est admissible de type fini. Donc il existe des caractres
${\chi}_1,..., \chi_n$ du centre
$Z$ de $M \cap \sigma(M)$ tels que pour tout $z
\in Z$, l'action de $(z - \chi_1(z))(z - \chi_2(z))... (z - \chi_n(z))$ 
soit nulle. En utilisant l'ŽgalitŽ:
$$ (z-\chi(z)b_z)(v \otimes b)= ((z-\chi(z))v)\otimes b, v \in V_{\delta}, b\in B$$
 on voit que: 
$(z -
\chi_1(z)b_z)...(z - \chi_n(z) b_z)$ agit trivialement sur le $(M \cap
\sigma(M), B)$-module $(H_0(M \cap \sigma(U), V_{\delta}))_{B}$, donc aussi
sur l'homologie lisse, $H_*(M \cap H, H_0(M \cap
\sigma(U),V_{\delta})_{B})$.\ste  Mais, d'aprs le Corollaire 2 (cf.  Proposition
\ref{actionG}) , pour tout cycle
$\varphi$,
$z_0\varphi$ est cohomologue ˆ $\varphi$ car $z_0$ est ŽlŽment du centre de $M
\cap H$. On en dŽduit que l'ŽlŽment $q$ de $B$ dŽfini par:$$q:= (1 -
\chi_1(z_0)b_{z_0})... (1 -
\chi_n(z_0) b_{z_0}) \in B$$ annule $H_*(M \cap H, H_0(M \cap \sigma(U),
V_\delta)_{B}$), donc $H_*(P \cap H, V_{\delta_B})$ d'aprs le Lemme 
prŽcŽdent. Alors $q$ a les propriŽtŽs voulues.\qed
\subsection{$(H,P)$-doubles classes} Un tore dŽployŽ de $G$,
$A_{\emptyset}$ contenu dans $\{ g \in G \vert \sigma(g) =
g^{-1}\}$ et de dimension maximale, sera dit tore $\sigma$-dŽployŽ
maximal,( ($\sigma,\F)$-torus dans [HWan]). On appelle $\sigma$-sous
groupe parabolique de $G$ tout sous-groupe, $P$,  de la forme $P_\lambda$, pour
un
$\lambda \in \Lambda( A_\emptyset)$, o $A_\emptyset$ est un tore
$\sigma$-dŽployŽ maximal. Alors  $\sigma(P) = P_{\lambda^{-1}}$ est
opposŽ ˆ $P$ relativement ˆ $A_{\emptyset}$.\ste  Enfin $HP$ est ouvert: en
effet les algbres de Lie de $P$ et $\sigma(P)$, ${ \p}$ et
$\sigma(\p)$ ont pour somme l'algbre de Lie $\g$ de $G$.\ste 
Soit $x \in \g$. On a  $x = y + z$ avec $y \in \p$, $z \in 
\sigma(\p)$. Donc $x = y' + h$ avec $y' = y - \sigma(z) \in \p$ 
et $h = z+\sigma(z) \in \h$. \ste Donc $\g = \h + \p$ 
ce qui implique que $HP$ contient un voisinage de $e$, l'ŽlŽment
neutre. Finalement $HP$ est bien ouvert.\ste D'aprs [HWan], Proposition 6.15, si
$P_0$ est un sous-groupe parabolique minimal de $G$, le nombre de
$(H,P_0)$-doubles classes est fini, donc aussi le nombre de
$(H, P)$-doubles classes.\ste Par ailleurs, d'aprs  [HH],
ThŽorme 2.9 (i), $P$ contient un $\sigma$-sous-groupe parabolique
minimal $P_{\emptyset}$ contenant $A_{\emptyset}$.
Le groupe $A_{\emptyset}$ est l'unique tore  $\sigma$-dŽployŽ
maximal du sous-groupe de LŽvi $\sigma$-stable de $P_{\emptyset}$,
$P_{\emptyset} \cap \sigma( P_{\emptyset})$.\ste
 On note $(A_i)_ {i \in I}$, un ensemble de reprŽsentants des classes de
$H$-conjugaison des tores $\sigma$-dŽployŽs maximaux de $G$.  On suppose
que cet ensemble contient $A_{\emptyset}$. Les $A_i$ sont tous conjuguŽs
sous
$G$, d'aprs la Proposition 1.16 de[HH].\ste On choisit, pour tout $i$, un
$x_i \in G$, avec $x_i A_{\emptyset} x_{i}^{-1} = A_i$ en prenant
$x_{\emptyset} = e$. On note ${\mathcal P}_i$ l'ensemble
des $\sigma$-sous-groupes paraboliques minimaux contenant $A_i$, qui est
fini (cf. [HH], Proposition 2.7).\ste  Les ŽlŽments de ${\mathcal P}_i$
sont tous conjuguŽs entre eux par un ŽlŽment du normalisateur de $A_i$ (
cf. l.c.). Comme les $A_i$ sont conjuguŽs entre eux, tous les ŽlŽments de
${\mathcal P}_i$ sont conjuguŽs sous $G$ ˆ $P_{\emptyset}$ et $P_i:= x_i$
$P_{\emptyset}$ $x_i^{-1}$.\ste 
On note $M_i$ le centralisateur dans $G$ de $A_i$. 
Si $L$ est un sous-groupe de $G$, on note $W_L(A_i)$ le quotient du
normalisateur dans $L$ de $A_i$ par son centralisateur. On note $W(A_i)$ au
lieu de 
$W_G(A_i)$.\ste 
On note ${\cal W}_i$, ou ${\cal W}_{i}^G$, un ensemble de reprŽsentants dans
$N_G(A_{\emptyset})$ de $W_{H_i}(A_{\emptyset})
\setminus W(A_{\emptyset})$ o $H_i = x_i^{-1} H x_i$. On note  ${\cal W}^G
= \bigcup_{i \in I} \{ x_i x \vert x \in {\cal W}_i^G\}$. Alors (cf.
[HH],  ThŽorme 3.1): 
\ber  \label{doubleclasseouverte}${\cal W}^G$ forme un ensemble de
reprŽsentants des $(H,P_{\emptyset})$-doubles classes ouvertes dans $G$. \eer
En particulier,  comme l'ensemble des $(H,P_{\emptyset})$-doubles classes est
fini (cf.  [HWan], Corollaire  6.16), on voit que 
$I$ est fini.\ste 
Soit $P$ un $\sigma$-sous-groupe parabolique contenant
$P_{\emptyset}$. On note $M$ le sous-groupe de Levi de $P$
contenant $A_{\emptyset}$.
On note ${\cal W}_{M,i}$ ou ${\cal W}_{M,i}^G$ un ensemble de
reprŽsentants dans $N_G(A_\emptyset)$ des doubles classes
$W_{H_i}(A_{\emptyset}) \setminus  W(A_\emptyset) / W_M(A_\emptyset)$ contenant
l'ŽlŽment neutre $e$.
 \begin{lem} \label{HPouverte} Toute $(H,P)$-double classe ouverte de $G$
est de la forme $HyP$ o $y$ est un ŽlŽment de:$${\cal W}_M^G = \bigcup_{i
\in I} \{x_i x \vert x \in{ \cal W}_{M,i}^G \}$$ En particulier toute
$(H,P)$-double classe ouverte est de la forme $HyP$ o $P^y$:= $yPy^{-1}$
est un
$\sigma$-sous-groupe parabolique de $G$.  On notera ${\overline {\cal W}}_M^G$
un ensemble de reprŽsentants des $(H,P)$-doubles classes ouvertes,  contenant $e$,  et
contenu dans 
${\cal W}_M^G$.
 \end{lem} 
\dem
 Soit $y = x_i x$ avec $x \in {\cal W}_i^M$. Alors
$Hx_ixP_\emptyset$ est ouvert d'aprs ce qui prŽcde. Donc $HyP$
est ouvert..\ste
 RŽciproquement si $\Omega$ est une $(H,P)$-double classe
ouverte, elle contient une $(H,P_\emptyset)$  double classe
ouverte (cf e.g. Appendice, Lemme \ref{orbites} ), donc de la forme $Hx_ixP$ avec $y \in
N_G(A_\emptyset)$ d'aprs (\ref{doubleclasseouverte}). Comme $Hx_iyP =
Hhx_iymP$ pour
$h
\in H$, $m \in M$, il en rŽsulte que $\Omega = H x_i xP$ pour un
$x \in {\cal W}_{M_i}$. \qed
Soit $P$ un sous-groupe parabolique $\sigma$-stable de $G$. On
note $A$ (resp $A_M$) le plus grand tore $\sigma$-dŽployŽ (resp. tore dŽployŽ) 
 du centre du sous-groupe de Levi $\sigma$-stable de $P$,
$M = P \cap \sigma (P)$. Si $x \in G$ et $E$ est une partie de $G$, on note
$E^x:= xEx^{-1}$.\ste  On note $X(M)_{\sigma}$ l'ensemble des caract\`eres non ramifi\'es de $M$ antiinvariants par $\sigma$. C'est un sous-groupe alg\'ebrique du groupe alg\'ebrique d\'efini sur $\C$, $X(M)$,  des caract\`eres non ramifi\'es de $M$.  On note $X$ (resp. $X_{\R}$) l'ensemble des caractres non
ramifiŽs de $M$ qui sont l'image par l'application de
(\ref{surjection}) (pour $M$ au lieu de $G$), de l'ensemble des
ŽlŽments de
$({\a_{M}^*})_{\C}$ (resp. $\a_{M}^*$) anti-invariants par
$\sigma$ ( qui laisse invariant  $A_M$  donc  aussi $ \a_M)$. Alors $X$ est  la composante connexe de l'\'el\'ement neutre du sous-groupe alg\'ebique $X(M)_{\sigma}$ de   $X(M)$. L'alg\`ebre $B_X$ est l'alg\`ebre des fonctions r\'eguli\`eres sur ce groupe alg\'ebrique connexe, qui est donc une vari\'et\'e irr\'eductible. Cette alg\`ebre est donc sans diviseur de z\'ero. On la note d\'esormais $B$. \ste 
On note $X_{\R}^x$ (resp.  $X^x$ ) l'ensemble des caractres
non ramifiŽs de $M^x$ obtenus par transport de structure de $M$ ˆ
$M^x$, via la conjugaison par $x$, des ŽlŽments de $X_\R$ (resp. $X$). \ste 
On note $A_{0}$ un tore dŽployŽ maximal $\sigma$-stable
contenu dans $P$ et contenant $A$.
\begin{lem} 
Avec les notations prŽcŽdentes, soit $\Omega$ une
$(H,P)$-double classe. Alors:\ste 
(i) $\Omega = HxP$ pour un $x$ tel que $A_{0}^x$ soit un tore
dŽployŽ maximal et $\sigma$-stable contenu dans $P^x$.\ste (ii)  Soit $\lambda$
un ŽlŽment de  $\Lambda(A^x)\subset\Lambda(A_{0}^x)$ tel
que $P^x = P_{\lambda}$. On note $\mu = \lambda
\sigma(\lambda)$. C'est un ŽlŽment du centre de $M^x \cap
\sigma(M^x)$. Supposons que  $\Omega$  ne soit pas ouverte. \ste 
Alors il existe $\chi \in X_{\R}^x$ tel que $\chi(\mu) \not= 1$. \end{lem}
\dem Ecrivons $\Omega = HyP$. D'aprs [HWan], Lemme 2.4, $P' := yPy^{ -1}$ contient
$A_{0}'$ un tore dŽployŽ maximal stable par $\sigma$. Soit $M'$ le
sous-groupe de LŽvi de $P'$ contenant $A_{0}'$. Alors $yMy^{-1}$ et
$M'$ sont des sous-groupes de LŽvi de $P'$, donc conjuguŽs par un ŽlŽment
$p'$ de $P'$. On peut mme se ramener au cas o $p'$ conjugue les tores
dŽployŽs maximaux de $P'$, $ yA_{0}y^{-1}$ et $A_{0}'$. 
On pose $x = p'y$. Alors $x A_{0} x^{-1} = A_{0}'$ et
$HxP = HyP$ car $p' \in yPy^{-1}$. 
Ceci prouve (i).\ste 
Montrons (ii). Comme $\Omega =HxP$ n' est pas ouvert, $HP^x$ n'est pas ouvert. Comme
$P^x=P_{\lambda}$  on n'a pas $\sigma(\lambda) = \lambda^{-1}$, sinon
$\sigma(P_{\lambda})$  serait un $\sigma$-sous-groupe parabolique de $G$ et
$HP_{\lambda}$ serait ouvert.
Donc $\mu$ est diffŽrent de $e$. Raisonnons par l'absurde et supposons
l'assertion du Lemme fausse.\ste Alors, par transport de structure, tout
ŽlŽment de $X$ est trivial sur $ x^{-1} \lambda x
x^{-1}\sigma(\lambda)x$.
  Mais $x^{-1} \lambda x \in \Lambda(A)$ et $x^{-1} \sigma(\lambda) x
\in \Lambda(A_{0})$. Le rŽseau $\Lambda(A_{0})$ s'identifie
gr$\hat{a}$ce ˆ l'application (\ref{H}) ( pour $G$ = $A_{0})$ ˆ un
sous-rŽseau dans $\a_{0}$.\ste  Notant
$\a_{\emptyset}$( resp. $\a_{0}^{\sigma}$) l'espace des ŽlŽments 
antiinvarants (resp. invariants)  de $\sigma$ dans
$\a_{
0}$, $\Lambda(A_{\emptyset})$ s'identifie ˆ un sous-rŽseau
de $\a_{\emptyset}$.\ste 
Gr\^ace ˆ (\ref{oplus}), $\a_M$ s'identifie ˆ un sous-espace de $\a_{0}$ et
$\a_{0} = \a_{0}^M \oplus \a_M$.\ste  Soit $\a:=  \a_{\emptyset} \cap  \a_M$ et soit 
$\a_M^{\sigma}$ l'espace  des points fixes de $\a_M$ sous
$\sigma$. Alors
$
\a_M = 
\a_M^{\sigma}
\oplus
\a$ et  $\Lambda(A)$ s'identifie ˆ un sous-rŽseau de $\a$.\ste  Il est
facile de voir qu'un ŽlŽment
$\Lambda( A_{0}) \subset
 \a_{0}$ est trivial sur tous les ŽlŽments de $X$ si et
seulement si sa projection sur $\a$ paralllement ˆ
$ \a_M^{\sigma} \oplus  \a_{0}^M$ est nulle. Notons
 $p_{\a}$ cette projection. Comme l'assertion de (ii) est supposŽe fausse, ce
qui prŽcde montre que:$$p_{\a}(x^{-1}
\lambda x ) = -p_{\a}(x^{-1} \sigma(\lambda) x )$$  et comme $x^{-1}
\lambda x \in \a$, on a:  \beq \label{pa} x^{-1} \lambda x = - p_{\a}(x^{-1}
\sigma(\lambda) x)\eeq On munit $ \a_{0}$ d'un produit scalaire
invariant par le groupe de Weyl $W(A_{0})$. On le prend mme invariant
par le groupe des automorphismes du systme de racines $\Sigma(A_{0})$
de $A_{0}$ dans l'algbre de Lie de $G$. 
En particulier le produit scalaire est $\sigma$-invariant. Alors la
projection $p_ {\a}$ est une projection orthogonale. Mais $x^{-1}
\lambda x$ et $\sigma(x^{-1}) \sigma(\lambda) \sigma(x)$ sont de mme
norme et:$$x^{-1} \sigma(\lambda) x = x^{-1} \sigma(x)(\sigma(x^{-1})
\sigma(\lambda) \sigma(x)) {\sigma(x)}^{-1} x$$
 Donc $x^{-1}
\sigma(\lambda) x$ et $\sigma(x^{-1}) \sigma(\lambda) \sigma(x)$ sont
conjuguŽs par $x^{-1} \sigma(x)$ qui normalise $A_{0}$ car $A_{0}^x$
est $\sigma$-invariant. Ces deux ŽlŽments de $ \a_{0}$ sont donc
de mme norme, Žgale ˆ celle de $x^{-1} \lambda x$, d'aprs les propriŽtŽs d'invariance
de la norme.\ste Comme
$p_ {\a}$ est une projection orthogonale, on  dŽduit de (\ref{pa}) que:  $$p_ {\a}(x^{-1}
\sigma(\lambda) x ) = x^{-1} \sigma(\lambda) x$$  et  $$x^{-1}
\lambda x = - x^{-1} \sigma(\lambda)x$$  dans $ \a_{0}$.\ste Ce
qui conduit ˆ:  $\lambda \sigma(\lambda) = e$. Une contradiction qui
achve de prouver le Lemme.\qed On conserve les notations du Lemme
prŽcŽdent. Soit $(\delta, V_{\delta})$ un reprŽsentation de $M$. On dŽfinit le $(M,B)$-module
$({\delta}_B, V_{\delta_B})$ comme en (\ref{deltaB}).\ste
On
notera $B^x$ l'algbre des fonctions sur $X^x$ dŽduite de $B$ par
transport de structure, via la conjugaison par $x$. Alors 
$((\delta_B)^x, V_{(\delta_B)^x}) $ a une
structure naturelle de
$(P^x,B^x)$-module Žquivalente par transport de structure ˆ
$((\delta^x)_{B^x}, V_ {(\delta^x)_{B^x}})$.
\begin{lem} Avec les notations ci-dessus, on suppose que $HxP$ n'est pas ouverte. \ste Alors
il existe
$q
\in B$, non constant, produit de facteurs de la forme $1-cb_m, c\in \C^*, m\in M$,  tel que la
multiplication par
$q$ annule le
$B$-module
$H_{*}(H, ind_{P^x
\cap H}^H {(\delta_B)^x}_{\vert P^x\cap H})$.\end{lem} 
\dem En appliquant le Lemme de Shapiro (cf. Proposition \ref{shapiro} ) et le
Lemme \ref{HcapP} on voit que le
$B$-module
$H_*(H, ind_{P^x \cap H}^H (\delta_B)^x)$ s'identifie ˆ
$H_*( M^x \cap \sigma(U), V_{^(\delta_B)^x})$.\ste Alors, par transport de
structure,  le Lemme rŽsulte de la Proposition \ref{z} appliquŽ ˆ $P^x$ et avec
$z_{0}=\mu$, car
$(\delta_B)^x$ s'identifie naturellement ˆ $(\delta^x)_{B^x}$. \qed
\subsection {$H$-homologie et $(H,P)$-doubles classes ouvertes} 
On note $\cal O$, la rŽunion des $(H,P)$-doubles classes ouvertes dans
$G$ et on dŽfinit:$$J_B = \{ \varphi \in I_B \vert supp \varphi \subset
\cal O\}$$ o $I_B = ind_P^G V_{\delta_B}$. C'est un sous-$(H,B)$-module de $I_B$.\ste
{\em On suppose dŽsormais que $\delta$ est admissible
de type fini.}
 
 \begin{theo}
\ste (i) Il existe $q\in B$, non nul, produit d'ŽlŽments  $1-c b_m, c\in \C^*, m\in M$,  
qui annule le
$B$-module
$H_{*}(H, I_B / J_B)$. \ste   (ii) On fixe un ensemble ${\overline{\cal W}}_M^G$ de
reprŽsentants des
$(H,P)$-doubles classes ouvertes comme dans le Lemme \ref{HPouverte}.
Alors:$$H_*(H, J_B) =
\bigoplus_{x \in {\overline{\cal W}}_M^G} H_*(M^x
\cap H, {\delta}^x) \otimes_{\C} B$$ En particulier $H_*(H,J_B)$ est un
$B$-module libre.\ste 
(iii)  Si $\xi \in Hom_B(H_0(H,J_B),B)$, il existe un unique
${\tilde \xi} \in Hom_B(H_0(H,I_B),B)$, tel que:$${\tilde \xi}(i(\psi)) = q
\xi(\psi),  \psi \in H_0(H,J_B)$$ o $i$ est le morphisme $H_0(H,J_B)
\to H_0(H,I_B)$, dŽduit de l'injection canonique $J_B \to I_B$. \end{theo}
\dem  (i) Le $H$-module $I_B/J_B$  admet, d'aprs la Proposition \ref{filtration}, une
filtration telle que l'on peut appliquer le Lemme 11 \`a ses  sous-quotients. La longue suite exacte
d'homologie permet de conclure que le produit des ŽlŽments de $B$,
dŽterminŽs par l'application de ce Lemme ˆ chacun des sous-quotients,
annule $H_*(H,I_B/J_B)$.\ste 
(ii) Notons, pour $x \in {\overline{\cal W}}_M^G$, $J_B^x = \{\varphi \in I_B \vert  supp
\varphi \subset HxP \}$, alors:$$J_B = \bigoplus_{x \in {\overline{\cal W}}_M^G}
J_B^x$$ Mais, d'aprs la Proposition \ref{filtration}, $J_B^x$ est isomorphe comme
$(H,B)$-module ˆ
$ind_{H\cap P^{x}}^{H}{(\delta_B)^x}_{\vert H\cap P^{x}}$.  Donc, d'aprs le Lemme de
Shapiro:$$H_*(H,J_B^x)
\simeq H_*(H
\cap M^x, V_{(\delta_B)^x})$$ 
et ceci comme  $B$-module. Mais d'aprs la dŽfinition de ${\overline{\cal W}}_M^G$ (cf.
Lemme \ref{HPouverte}), si
 $x \in {\overline{\cal W}}_M^G$:$$xAx^{-1} \subset \{ g \in G  {\vert}  \sigma(g) =
g^{-1} \}$$ et c'est un tore $\sigma$-dŽployŽ maximal du centre de
$M^x$.  Alors si $x \in {\overline{\cal W}}_M^G$, $P^x$ est
un $\sigma$-sous-groupe parabolique et $M^x$ est $\sigma$-stable, donc Žgal
au sous-groupe de LŽvi $\sigma$-stable de $P^x$. Les ŽlŽments de $X$ se
transportent par $x$ en des caractres non ramifiŽs de $M^x$, triviaux
sur $H \cap M^x$. Il en rŽsulte que comme $(H \cap M^x, B )$-module,
$V_{(\delta_B)^x}$ est isomorphe ˆ $V_{\delta^x} \otimes B$ o $H \cap
M^x$ n'agit que sur le premier facteur et $B$ que sur le second. (ii)
en rŽsulte immŽdiatement.\ste
Prouvons (iii). On dispose de la suite exacte de $B$-modules: \beq
\label{suitecourte}  H_1(H,I_B/J_B) \opto{c} H_0(H,J_B) \opto{i} H_0(H,I_B)
\opto{p} H_0(H,I_B/J_B) \to 0  \eeq qui est donnŽe par la longue suite exacte
d'homologie. Ici $i$ est donnŽ par passage au quotient de l'injection
naturelle $J_B \to I_B$, et $p$ par passage au quotient de la
projection naturelle $I_B \to I_B/J_B $.\ste
Soit $\varphi \in H_0(H,I_B)$. Alors $q\varphi$ a une image nulle
dans $H_0(H,I_B/J_B)$, d'aprs (i). Donc $q\varphi$ est dans l'image par $i$
de $H_0(H,J_B)$, d'aprs (\ref{suitecourte}),  i.e.:\beq   q\varphi =
i(\psi)\hbox{  pour un }
\psi
\in H_0(H,J_B)\eeq   
Montrons que $\psi$ est unique. En effet si:  $i(\psi) = i(\psi') =
q\varphi$, alors $\psi - \psi'$ est dans le noyau de $i$, donc dans
l'image de $c$, i.e.:   $\psi - \psi' = c(\eta)$ pour $\eta \in
H_1(H,I_B/J_B)$. D'aprs (i), $\eta$ est annulŽ par la multiplication
par $q$. Donc $q\eta $ est nul  ainsi que $c(q\eta) = q(\psi - \psi')
\in H_0(H,J_B)$. Mais d'aprs (ii), $H_0(H,J_B)$ est un $B$-module
libre. Donc $\psi - \psi' = 0$ comme dŽsirŽ.\ste 
On pose alors:\ber ${\tilde \xi} (\varphi):= \xi(\psi)$  pour $\varphi \in
H_0(H,I_B)$,  o $\psi$ est l'unique ŽlŽment de $H_0(H,J_B)$ tel que 
$i(\psi) = q\varphi$\eer 
Calculons maintenant ${\tilde \xi}(i(\psi))$ pour $\psi \in H_0(H,J_B)$.
Comme   $qi(\psi) = i(q\psi)$, on a:$${\tilde \xi}(i(\psi)) =
\xi(q\psi) = q\xi(\psi)$$comme dŽsirŽ.\ste Si un autre ${\tilde \xi}'$
vŽrifiait (iii), alors ${\tilde \xi} - {\tilde \xi}'$ serait nul sur
$i(H_0(H,I_B))$ qui est Žgal ˆ $Ker(p)$. D'aprs (i), ce noyau contient
$qH_0(H,I_B)$. Donc:$$({\tilde \xi} - {\tilde \xi}')(q\varphi) =
q({\tilde \xi} - {\tilde \xi}')(\varphi) = 0, \>\> \varphi \in
H_0(H,I_B)$$ Comme l'algbre $B$ est  sans diviseur de 0 ( voir apr\`es le Lemme 9),  
$({\tilde \xi} - {\tilde \xi}')(\varphi) \in B$ est nul  pour tout $\varphi \in
H_0(H,I_B)$   comme dŽsirŽ.\qed
\subsection {Localisation du ThŽorme 1 }  Soit $K$ un sous-groupe
compact maximal de $G$ qui est le fixateur d'un point spŽcial de
l'appartement associŽ ˆ $A_{0}$ dans l'immeuble de $G$.\ste 
On note $I$ l'espace de  $ind_{K \cap P}^K \delta _{\vert_{K
\cap P}}$. Alors la restriction des fonctions ˆ $K$ dŽtermine un
isomorphisme de $K$-modules entre l'espace $I_{\chi}$ de  $\pi_{\chi}:=ind_P^G (\delta \otimes
\chi )$, $\chi \in X(M)$   (resp. $I_B$ de $\pi_{B}:=ind_P^G \delta_B$)   et $I$
(resp.  $I \otimes B$). On note $ \overline {\pi}_{\chi} $ (resp $ \overline
{\pi}_{B} )$ la reprŽsentation de $G$ sur $I$ (resp.  $ I \otimes B$) 
 dŽduite de $\pi_{\chi}$ ( resp $\pi_{B}$ ) par transport de structure via cet
isomorphisme.\ste 
Si $\varphi \in I$ et $\chi \in X(M)$, on note $\varphi_{\chi}$ l'ŽlŽment de
l'espace $I_{\chi}$ correspondant ˆ $\varphi$ par cet
isomorphisme.\ste  On note $\varphi_{B}$ l'ŽlŽment de l'espace $I_B$ de $ind_P^G \delta_B$ dont
la restriction ˆ $K$ est $\varphi \otimes 1_B \in I \otimes B$, o $1_B$ est
l'unitŽ de $B$. On note
$ev_{\chi}$ l'Žvaluation en $\chi$ des ŽlŽments de $B$  et $\cal J_\chi$ son noyau. Alors
$$\varphi_\chi = ev_\chi(\varphi_B)$$ o on a notŽ $ev_\chi(\varphi_B)$ au lieu de
$(id_{V_\delta} \otimes ev_\chi) ( \varphi_B) $. Il en rŽsulte que: \ber Pour tout
$\varphi \in I$ et $g \in G$, $\chi\mapsto {\overline\pi}_{\chi}(g) \varphi$ est un
ŽlŽment de
$I$ $\otimes$ $B$  Žgal ˆ ${\overline \pi}_B(g)(\varphi \otimes 1_B) $.\eer 
Pour tout $B$-module $M$ on note $M_{(\chi)}$ = $M/{\cal J}_\chi M$ sur lequel tout ŽlŽment
$b$ $\in$ $B$ agit par $b(\chi)$. Alors on a facilement: \ber l'application $ev_{\chi}$
entre  $I_B= ind_P^G V_{\delta _B}$ et $ind_P^G V_{\delta \otimes \chi}$ passe au quotient en
un isomorphisme entre  
$(ind_P^G V_{\delta _B})_{(\chi)}$ et $ind_P^G V_{\delta \otimes \chi}$,  \eer
 ce qui se voit
aisŽment dans la rŽalisation compacte.

On note $J_\chi = \{ \varphi \in I_\chi \vert Supp \varphi \in \cal O\}$
\begin{theo} (i) Pour tout $\chi \in X, n \in \N$, on a un isomorphisme linŽaire
naturel:$$H_n(H,J_\chi) \simeq \oplus_{x \in {\overline{\cal W}}_M^G} H_n(M^x \cap H,
V_{\delta^x})$$ (ii)  Faisant $n = 0$ dans (i) et passant au dual, on en dŽduit un
isomorphisme naturel: $$ J_\chi^{*H} \simeq {\cal V}_{\delta} \hbox{ o }{\cal V}_{\delta}:
=  
\oplus _ {x \in {\overline{\cal W}}_M^G} ({V_{\delta^x}^*})^{M^x \cap H}$$ 
Celui-ci  associe ˆ $\eta \in  {V_{\delta^x}^*}^{M^x \cap H}$  (i.e  $\eta \in
{V_\delta^*}^{M \cap x^{-1}Hx}$), $\xi(P,\delta,\chi,\eta) \in J_\chi^*$ dŽfini par
$$\xi(P,\delta,\chi,\eta)(\varphi) = \int_{H/M^x \cap H} <\varphi(hx),\eta> d\dot{h},
  {\varphi \in J_\chi}$$ \ste  
(iii)  Avec les notations du ThŽorme 1, soit $\chi \in X$ avec $q(\chi) \ne 0$.
Notons $i_\chi$ l'injection naturelle de $ J_\chi $ dans $ I_\chi$. Celle-ci
passe au quotient en un isomorphisme notŽ encore $i_\chi: $$H_0(H,J_\chi) \tilde{\to}
H_0(H,I_\chi)$ dont la transposŽe dŽtermine un isomorphisme: ${I^*_\chi}^H \to
{J^*_\chi}^H$ donnŽ par la restriction, $r_{\chi}$,  des formes linŽaires de $I_\chi$ ˆ
$J_\chi$.\ste On notera pour $\eta \in {\cal V}_\delta $, $j(P,\delta,\chi,\eta)$ l'ŽlŽment
de
${I^*_\chi}^H$ correspondant ˆ $\xi(P,\delta,\chi,\eta)$ dans cet isomorphisme.\ste 
(iv) On note ${\overline j}(P,\delta,\chi,\eta)$ l'ŽlŽment de $I^*$ dŽduit de
$j(P,\delta,\chi,\eta)$ par transport de structure ˆ l'aide de $r_\chi$. Alors, pour
tout $\varphi \in I$, l'application: $\chi \mapsto
q(\chi)<\bar{j}(P,\delta,\chi,\eta), \varphi>$ dŽfinie sur l'ensemble des $\chi \in
X$ tels que $q(\chi) \not= 0$ admet un prolongement rŽgulier ˆ $X$ tout entier (i.e.
se prolonge en un ŽlŽment de $B$).\ste On dira que $\chi \mapsto q \bar{j}
(P,\delta,\chi,\eta)$ est une famille polynomiale de formes linŽaires sur $I_\delta$.    
\end{theo}
 \dem 
(i) se prouve comme la partie (ii) du ThŽorme 1. \ste 
(ii) rŽsulte de l'explicitation des isomorphismes. \ste
Prouvons (iii). Montrons la surjectivitŽ de $r_\chi$. Soit $\theta \in J_\chi^{*H}$.
D'aprs l'isomorphisme du ThŽorme 1 (ii), il existe $\xi \in Hom_B(H_0(H,J_B), B)$ tel que
$\xi$ passe au quotient en un ŽlŽment de: $$Hom_\C(H_0(H,J_B)_{(\chi)}, B_{(\chi)}) \simeq
Hom_\C(H_0(H,J_\chi),\C) \simeq J_{\chi}^{*H}$$
Žgal ˆ $\theta$. Alors utilisant l'ŽlŽment $\tilde{\xi}$ de $Hom_B(H_0(H,I_B), B)$ donnŽ
par le ThŽorme 1 (iii) et localisant en $\chi$ comme ci-dessus, on trouve $\tilde{\theta}
\in I_{\chi}^{*H}$ tel que:  $$\tilde{\theta}(\psi) = q(\chi) \theta(\psi), \>\>
\psi
\in J_\chi$$ 

Alors si $ q(\chi) \not= 0$, $q(\chi)^{-1} \tilde{\theta} \in  I_\chi^{*H}$  a pour
image $\theta$ par $r_\chi$. Donc $r_\chi$ est bien surjective.\ste 
Montrons maintenant que $r_\chi$ est injective. Soit ${\tilde \theta} \in Ker\space \>
r_\chi$. Alors ${\tilde \theta}$ s'identifie ˆ un ŽlŽment de: 
$$(I_\chi/J_\chi)^{*H} \simeq H_0(H,I_\chi/J_\chi)^*$$
Comme $(I_B/{\cal{J}_\chi}I_B)/(J_B/{\cal{J}_\chi}J_B) \simeq I_B/(J_B +
{\cal{J}_\chi}I_B)$, on a un isomorphisme de $H$-modules: $$
 I_{\chi}/J_{\chi}
\simeq (I_B/J_B)_{(\chi)}$$
Donc  $H_0(H, I_\chi/J_\chi) \simeq H_0(H,{(I_B/J_B)}_{(\chi)})$. Mais on vŽrifie facilement
que pour tout $(H,B)$-module lisse, $V$, on a : 
\beq  \label{I/Jchi} (H_{0}(H, V))_{(\chi )}\simeq H_{0}(H,V_{(\chi)})\eeq

Donc 
$$H_0(H,I_\chi/J_\chi) = H_0(H,I_B/J_B)_{(\chi)}$$  
 Mais la multiplication par $q$ annule
$I_B/J_B$ donc aussi $H_0(H,I_B/J_B)$. Comme $q(\chi) \not= 0$, on en dŽduit que 
$H_0 (H,I_{B}/J_{B})_{(\chi)} = \{0\}$ et l'isomorphisme prŽcŽdent implique que 
$H_0(H,I_\chi/J_\chi)^*$ est rŽduit ˆ $\{0\}$. Donc ${\tilde \theta} = 0$ et $r_\chi$
est bien injective.\ste 
Prouvons (iv). Soit $\eta \in \cal{V}_\delta$. D'aprs le thŽorme 1 (ii), 
$Hom_B(H_0(H,J_B), B)$ s'identifie naturellement ˆ ${\cal V}_\delta$. 
On note $\xi$ l'ŽlŽment de $Hom_B(H_0(H,J_B), B)$ correspondant ˆ
$\eta$ par  cet isomorphisme et $\tilde{\xi}$  l'ŽlŽment de
$Hom_B(H_0(H,I_B), B)$ dŽterminŽ par le ThŽorme 1 (iii). \ste 
L'ŽlŽment $\tilde{\xi}$ de $Hom_B(H_0(H,I_B ),B)$ dŽtermine, ˆ l'aide de   la
projection naturelle $I_B \to H_0(H,I_B)$ et de l'isomorphisme de $B$-modules,
dŽcrit plus haut, entre $I_B$ et
$I\otimes B$, un ŽlŽment $\bar{\xi}$ de $Hom_{B}(I_\delta \otimes B,
B)$.\ste Les constructions (cf. la preuve de la surjectivitŽ de $r_{\chi}$ dans (iii))
montrent que, lorsque $q(\chi)$ est non
nul:$$<q(\chi)\bar{j}(P,\delta,\chi,\eta), \varphi> = 
(\bar{\xi}(\varphi \otimes 1_B))(\chi)$$
Alors (iv) rŽsulte du fait que  $\bar{\xi}( \varphi \otimes 1_B) \in
B$.  \qed 
\subsection {ReprŽsentations rationnelles $H$-sphŽriques} Les deux faits suivants montrent  
 comment des questions d'invariance par un sous-groupe se ramnent ˆ des questions
d'invariance par une sous-algbre de Lie.     On note  ${\underline G}$ un   groupe algŽbrique 
dŽfini sur un corps algŽbriquement clos
$\k$, de caractŽristique zŽro.  Alors (cf. [Hu2], 13.1): 
\ste Soit $(\pi, V)$ une reprŽsentation rationnelle de dimension finie de ${\underline G}$ et
$W$ un sous-espace vectoriel de $V$. On note encore $\pi$ la reprŽsentation de l'algbre de Lie
$\g$ de ${\underline G}$ dans $V$. Alors le sous-groupe
$H:=\{g\in  {\underline G} \vert \pi(g)W\subset W\}$ est un sous-groupe fermŽ de ${\underline
G}$, admettant
$\h:=\{X\in
\gÊ\vert \pi(X)  W\subset W\}$ comme algbre de Lie. 
\ste
De plus (cf. [Hu2], 13.2)):  
\ste  Si $H$ est un sous groupe fermŽ de ${\underline G}$, possŽdant la mme algbre de Lie que
${\underline G}$, alors 
$H={\underline G}$.\ste 
Nous utiliserons ces rŽsultats pour la cl™ture algŽbrique de $\F$. Nous utiliserons
Žgalement les  propriŽtŽs des reprŽsentations de plus haut poids de dimension finie pour les
algbres de Lie semi-simples  sur cette cl™ture algŽbrique (cf. [Hu1], ch. VI). \ste\ste    
Soit
$P_{\emptyset}$ un
$\sigma$-sous-groupe parabolique minimal de $G$, $M_{\emptyset}$ = $P_{\emptyset}\cap$
$\sigma(P_{\emptyset})$ le sous-groupe de LŽvi $\sigma$-stable de
$P_{\emptyset}$,
$A_{\emptyset}$ le plus grand tore $\sigma$-dŽployŽ du centre de $M_{\emptyset}$ qui est
un tore $\sigma$-dŽployŽ maximal. On fixe $A_{0}$ un tore dŽployŽ maximal
contenant
$A_{\emptyset}$, i.e. un tore dŽployŽ maximal de $M_{\emptyset}$  et
$\sigma$-stable.\ste  On remarque qu'alors $\sigma$ agit naturellement
sur 
$\a_{0}$  et que
$\a_{\emptyset}$  s'identifie au sous-espace des ŽlŽments antiinvariants de $\a_{0}$. On
note $P_{0}$ un sous-groupe parabolique minimal de $G$ contenu dans
$P_{\emptyset}$ et contenant $A_{0}$. 
On note $P$ un $\sigma$-sous-groupe parabolique de $G$
contenant $P_{\emptyset}$ et $M$ son sous-groupe de LŽvi $\sigma$-stable. On note
$U_{0}$ (resp. $U_{\emptyset}$, resp. $U$) le radical unipotent de $P_{0}$ ( resp.
$P_{\emptyset}$, resp.  $P$ ). On note $A_{G,\sigma}$ le plus
 grand tore $\sigma$-dŽployŽ de
$A_G$ et  ${\a}_G^{\sigma}$ (resp. ${\a}_{G,\sigma}$ )
l'ensemble des points fixes ( resp. antiinvariants ) de ${\a}_G$ sous
$\sigma$. \ste 
On note $p_{\sigma}$ la projection de ${\a}_G$ sur $
\a_{G,\sigma}$ paralllement ˆ ${\a}_G^{\sigma}$ et $H_{G,\sigma}$
la composŽe $p_{\sigma} \circ H_G$.\ste 
Alors, remarquant que tout caractre rationnel  de $G$
tel que $\chi^{\sigma} =
\chi^{-1}$ vŽrifie $\chi^2= \chi \chi^{-\sigma}$,   on Žtablit immŽdiatement
que:\ber Le noyau
$G_{\sigma}^1$ de
$H_{G,\sigma}$ est Žgal ˆ l'intersection des noyaux des caractres ${\vert
\chi \vert}_{\F}$ o $\chi$ dŽcrit l'ensemble des caractres rationnels de $G$
tels que $\chi^{\sigma} =
\chi^{-1}$.\eer  \ste 
Nous appelerons dans cet article reprŽsentation de plus haut poids $\Lambda$, 
$\Lambda
\in Rat(M_{0})$,  une reprŽsentation rationnelle de
$G$, dŽfinie sur $\F$   dans un espace vectoriel sur $\F$ de dimension finie, ({\em non
nŽcessairement irrŽducible}) possŽdant un vecteur non nul $v_{\Lambda}$, dit de plus
haut poids
$\Lambda$, de poids $\Lambda$ et invariant par
$U_{0}$ et se transformant par
$\Lambda$ sous $M_{0}$, et telle que les poids de $A_0$ dans cette reprŽsentation,
identifiŽs ˆ des ŽlŽments de $\a_0^*$ soient Žgaux  ˆ la diffŽrence de la forme linŽaire
correspondant ˆ
$\Lambda$ avec une somme ˆ coefficients entiers positifs, de racines de $A_0$ dans
l'algbre de Lie de
$U_0$.\ste Par passage ˆ la cl™ture algŽbrique et ˆ l'algbre de Lie (voir le dŽbut du
paragraphe)  on voit qu'une
 reprŽsentation irrŽductible de plus haut poids $\Lambda$, si elle existe, est unique ˆ
isomorphisme prs. De mme on voit qu'une
reprŽsentation rationnelle de dimension finie est irrŽductible de plus haut poids
$\Lambda$, si et seulement si elle est engendrŽe par un vecteur de plus haut poids
$\Lambda$.
\ber
\label{caractererationel} On remarque, par passage ˆ l'algbre de Lie,  que si
$\Lambda
\in Rat(M_{0})$ est le plus haut poids d'une reprŽsentation irrŽductible  rationnelle de
$G$ et
$\Lambda \in {\a}_G^*$, alors cette reprŽsentation  est un caractre rationnel de $G$.\ste 
Si, de plus, $\Lambda \in {\a }_{G,\sigma}^*$, ce caractre,  $\chi_{\Lambda}$,  vŽrifie
$\chi_{\Lambda}\circ \sigma= \chi_{\Lambda}^{-1}$.   \eer Soit $\Lambda\in
Rat(M_0)$ {\em tel que
$\Lambda
\in
\a_M^{*}$, i.e. tel qu'il soit  nul sur $\a^M_{M_0}$}. \ste Pour une  reprŽsentation
irrŽductible  de plus haut poids
$\Lambda$ ( si elle existe), $(\pi_{\Lambda}, V_{\Lambda})$, on note $$e_{\Lambda,H}^* \in (
V_{\Lambda} \otimes V_{\Lambda}^*)^*
\simeq V_{\Lambda}^* \otimes V_{\Lambda} \simeq {End V}_{\Lambda} $$ l'ŽlŽment
correspondant ˆ l'identitŽ de $V_{\Lambda}$.\ste 
Alors $e_{\Lambda,H}^* = \sum_{i = 1}^m e_i^* \otimes e_i $ o $(e_i)$ est une base de
$V_{\Lambda}$ et $(e_i^*)$ la base duale. On note $v_{\Lambda}^*$ un ŽlŽment de
$V_{\Lambda}^*$ de poids ${\Lambda}^{-1}$ sous $M_{0}$ et vŽrifiant $<v_{\Lambda}^*,
v_{\Lambda} > = 1$. Celui-ci est unique. La droite $\F v_{\Lambda}^*$ est
invariante par le sous-groupe parabolique opposŽ ˆ $P_{0}$, $P_{0}^{opp}$, relativement
ˆ
$A_0$.\ste Alors, d'aprs (\ref{caractererationel}) appliquŽ ˆ $M$ au lieu de $G$, on
voit que 
$\F v_{\Lambda}^*$ est invariante par $M$, car $\Lambda\in \a_M^*$ . Donc cette
droite est aussi invariante par le sous-groupe engendrŽ par $M$ et $P_{0}^{opp}$, qui
est Žgal au sous-groupe parabolique opposŽ ˆ $P$ , relativement ˆ $A_0$. Mais celui-ci
est Žgal ˆ
$\sigma(P)$ puisque $P$ est un
$\sigma$-sous-groupe parabolique contenant $P_{\emptyset}$. Alors $v_{\Lambda}^*$ est un
vecteur de plus haut  poids $\Lambda^{-1}\circ \sigma$ pour la reprŽsentation
irrŽductible
$\pi_{\Lambda}^{*}\circ \sigma$. 
   En
prenant
$e_1 = v_{\Lambda}$ et les autres
$e_i$ dans le noyau de
$v_{\Lambda}^*$, on a $e_1^* = v_{\Lambda}^*$ et on voit que: $$< e_{\Lambda,H}^*,
\tilde{v}_{\Lambda} > = 1 , \>\> avec \>\> \tilde{v}_{\Lambda} = v_{\Lambda} \otimes
v_{\Lambda}^*$$ Par ailleurs
$e_{\Lambda,H}^*$ est invariant sous $G$ par ${\pi}_{\Lambda}^* \otimes
{\pi}_{\Lambda}$, donc par $H$  sous ${\pi}_{\Lambda}^* \otimes ({\pi}_{\Lambda} \circ
\sigma)$. On notera:  
$$\tilde{\Lambda}:=\Lambda
(\Lambda^{-1}\circ
\sigma), ({\pi}_{\tilde{\Lambda}},V_{\tilde{\Lambda}}):=({\pi}_{\Lambda}
\otimes ({\pi}_{\Lambda}^*\circ \sigma ),  V_{\Lambda} \otimes V_{\Lambda}^*) $$

Une inspection
des poids montre que l'espace de poids
$\tilde{\Lambda}$ de $V_{\tilde{\Lambda}}$  sous $M_{0}$ est de dimension $1$, et 
la reprŽsentation ${\pi}_{\tilde{\Lambda}}$,{\em qui n'est pas
nŽcessairement irrŽductible}, est de plus haut poids $\tilde{\Lambda}$. 
\ste  
 \ber $(\pi_{\tilde{\Lambda}} $,$ V_{{\tilde\Lambda}}$) est une reprŽsentation de plus
haut poids ${\tilde\Lambda}$, avec un vecteur de plus haut poids ${\tilde\Lambda}$,
$v_{{\tilde\Lambda}}$ et un vecteur $H$-invariant dans
$({ V}_{{\tilde \Lambda}})^*$ pour $({\pi_{\tilde\Lambda}})^*$,
 ${\tilde e}_{ \Lambda,H}^*$ vŽrifiant $<{\tilde e}_{ \Lambda,H}^*,v_{\tilde \Lambda} >
= 1$
\eer 
Un ŽlŽment de $ \a_{0}^* $ antiinvariant par $\sigma$ est un  ŽlŽment de
${\a_\emptyset}^*$ par la maximalitŽ de $\a_{\emptyset}$.\ste 
On note $\Sigma(G, A_{0})$ (resp. $\Sigma(P_{0}, A_{0})$ ou $\Sigma(P_0)$) l'ensemble des
racines   de
$A_{0}$ dans l'algbre de Lie de $G$, (resp. $P$). On note $\Delta(P_{0})$ l'ensemble des
racines simples de $\Sigma(P_{0})$. Si $\Theta$ est une partie de $\Delta(P_{0})$,
on note $<\Theta>$ le sous-systme de $\Sigma(P_{0})$ engendrŽ par $\Theta$ et
$P_{<\Theta>}$ le sous-groupe parabolique de $P_{0}$ pour lequel $\Sigma(P_{0}) \cup
<\Theta>$ est l'ensemble des racines de $A_{0}$ dans l'algbre de Lie de $P_{\Theta}$.
On dŽfinit   $\Theta_{\emptyset}, \Theta\subset \Delta(P_{0}) $   par les ŽgalitŽs:
$$P_{\emptyset}= P_{<\Theta_{\emptyset}>}, P=P_{<\Theta>}$$ On Žcrit  
$$ \Delta (P_0)=\{\alpha_1, \dots, \alpha_k\},\Delta (P_0)\setminus \Theta_{\emptyset } =  \{\alpha_1, \dots,
\alpha_l\}, 
 \Delta (P_0)\setminus \Theta =  \{\alpha_1, \dots, \alpha_m\},
   $$
de sorte que $k\geq l \geq m$. 
\ste On note $\delta_1,...,\delta_k \in \a_{0}^*$ les poids fondamentaux. Ils  sont nuls
sur
$\a_G$, et pour $i=1,\dots, m$, $\delta_i\in \a_M^*$.
\begin{prop}  \label{pitilde}\ste (i)  Il existe des entiers ${n_1,...,n_k} \in {\N}^*$
tels que
$n_i\delta_i$ corresponde ˆ un plus haut poids $\Lambda_i$ d'une reprŽsentation rationnelle
$(\pi_i, V_i)$ de $G$. \ste (ii)  Avec les notations prŽcŽdentes, on note ${\tilde
V}_i$ au lieu de
$V_{\tilde{\Lambda}_i}$, $\tilde{\pi_i}$ au lieu de $\pi_{\tilde{\Lambda}_i}$,...\ste 
Pour $i = 1,...,m$, la droite $\F\tilde{v}_i$ est invariante par $M$ et $\tilde{v}_i$ est
invariant par $M_\sigma^1$.
\end{prop}  
On notera $${\tilde \delta}_i:= \delta_{i}-\delta_{i}\circ \sigma
$$
\dem  La partie (i) de la
proposition rŽsulte de la thŽorie des reprŽsentations rationnelles de
$G$ (cf. [BoT],proposition 12.13).  Passons ˆ (ii). Soit
$i
\in
\{1,...,m\}$. Alors
$\tilde{v}_i$ engendre une reprŽsentation de $M$ de plus haut poids
$n_i\tilde{\delta_i}$ qui est un ŽlŽment de $\a_M^*$ car $i \in
\{1,...,m\}$ et mme de $\a_{M,\sigma}^*$,   auquel on peut appliquer
(\ref{caractererationel}).
\qed               
\begin{rem} (i)  Un argument utilisant les puissances tensorielles montre que les multiples
entiers non nuls des $n_i$ vŽrifient les mmes propriŽtŽs.\ste 
(ii) Ecrivant $\a_0 = \a_0^\sigma \oplus {\a}_{\emptyset}$, on a $\tilde{\delta}_i \in
\a_\emptyset^*$ \space pour $i = 1,...,l$. De plus Žcrivant $\a_\emptyset = (\a_0^M
\cap  \a_\emptyset) \oplus \a_{M,\sigma}$, on a $\tilde{\delta}_i \in
\a_{M,\sigma}^*$ \space pour $i =1,...,m$ \space avec $\tilde{\delta}_i$ nul sur
$\a_{G,\sigma}$. De plus l'ensemble $\{\tilde{\delta}_i \vert i = 1,...,m\}$ engendre
$(\a_{M,\sigma}/\a_{G,\sigma})^*$. On suppose la numŽrotation choisie de sorte que
$(\tilde{\delta}_1,...,\tilde{\delta}_{m'})$ forme une base de $\a_{M,\sigma}^*$ avec $m'
\le m$.\end{rem} \begin{prop}\label{deltaM}
  (i) Soit  $i \in \{m+1,...,l\}$, i.e. $\alpha_{i}\in \Theta\setminus
\Theta_{\emptyset}$. Il existe un entier
$n_i$ tel que
$n_i\delta_i$ (resp ${n_i}{\delta_i}^{M}$,  qui est l'ŽlŽment de $\a_0^{*}$   Žgal
ˆ $n_{i}\delta_i $ sur $\a_{0}^{M}$ et nul sur $\a_{M}$),  corresponde ˆ un plus haut
poids
$\Lambda_i$ (resp $\Lambda_i^M$) d'une reprŽsentation rationelle de $G$ (resp. de $M$).
Les multiples entiers des $n_i$ vŽrifient les mmes propriŽtŽs. Noter que $\delta_i^{M}$ 
est
un poids fondamental pour le systme de racines $<\Theta >$ et $(\Lambda_i)
(\Lambda_i^M)^{-1}$ dŽtermine, gr\^ace ˆ (2.33), un caractre rationnel $\Lambda_{i,M}$
de
$M$, qui  correspond ˆ
$n_i{\delta_i}_{\vert \a_{M}} \in \a_{M}^*$. 
\ste (ii) Le vecteur
$v_{\Lambda_i}$ engendre sous $M$ une reprŽsentation irrŽductible de plus haut
poids $ \Lambda_i$  notŽe ${\overline V}_i^M$ contenue 
dans $V_{i}$. Alors ${\tilde {\overline  V}}^M_i:={\overline V}_i^M\otimes({\overline
V}_i^M)^*$ s'idenditifie ˆ un sous espace de ${\tilde V}_i$, qui est stable par ${{\tilde
\pi}_i}{}_{\vert M} $. Cette reprŽsentation s'identifie au produit tensoriel du caractre
rationnel
${\tilde
\Lambda_{i,M}}:=
\Lambda_{i,M} ( 
\Lambda_{i,M}^{-1}\circ
\sigma)$ de
$M$ avec la reprŽsentation irrŽductible de
$M$,
${\tilde V}_i^M$, associŽe ˆ  $ \Lambda_i^M$ par la construction (2.24).
\ste La restriction de ${\tilde e}_{i,H}^*$ ˆ ${\tilde {\overline  V}}^M_i$
dŽtermine un vecteur non nul,  $H\cap M$-invariant,  qui s'identifie au
vecteur
$1\otimes {\tilde e}_{i,M\cap H}^*$.
\end{prop}
\dem (i) rŽsulte de la Proposition et de la Remarque prŽcŽdentes. \ste 
(ii) est immŽdiat.\qed
\subsection{PropriŽtŽs de l'adhŽrence des $(H,P)$-doubles classes ouvertes}
On choisit les entiers $n_i$ comme dans les Propositions prŽcŽdentes.\ste 
Pour $i =1,...,m$, on dŽfinit: $$\varepsilon_i(g) = \vert <{\tilde
\pi}_i(g){\tilde v}_i,{\tilde  e}_{i, H}^*>\vert_{\F},  g
\in G$$ On note $\varepsilon = \prod_{i=1}^{m'} \varepsilon_i$.
\begin{prop} (i) L'ouvert de $G$, $HP$, est contenu dans $\{ g \in G \vert  \space
\varepsilon(g) \not= 0 \}$. \ste 
(ii) Soit $\overline{HP}$ l'adhŽrence de $HP$ dans $G$. Alors $\overline{HP} \backslash HP$
est contenu dans $\{ g \in G \vert  \varepsilon(g) = 0\}$.
\ste(iii) Si la suite $(g_n) = (h_nm_nu_n)$, avec $h_n \in H, m_n \in M, u_n
\in U$, converge vers un ŽlŽment de $\overline{HP} \backslash HP$ et $\nu \in
\a^*_{M,\sigma}$ est strictement
$P$-dominant, alors la suite
$(e^{      \nu(   H_{M,\sigma}(m_n) )     }) $ tend vers zŽro.\end{prop}
\dem Elle est analogue ˆ celle du Lemme 4 de [BrD]. Nous allons
expliquer les modifications nŽcŽssaires.\ste 
D'abord si $g=hmu $ avec $h\in H$, $m\in M$, $u\in U$, $\varepsilon(g)$ est Žgal ˆ:
$$ \prod_{i=1}^{m'}\vert < {\tilde  \pi}_i(m) {\tilde v}_{i}, {\tilde e}_{i, H}^*
>\vert_{\F}$$ qui est non nul  d'aprs la Proposition  \ref{pitilde} (ii). Ceci achve de
prouver (i). 
\ste Montrons (ii). On considre une suite $(g_n)$ dans $HP$ convergeant vers $g \in
\overline{HP} \setminus  HP$.  On Žcrit $g_n = h_nm_nu_n$ avec $h_n \in H, m_n \in M, u_n
\in U$. On va montrer que l'image
dans $\a_{M,\sigma}/ \a_{G,\sigma}$ de la suite $(H_{M,\sigma}(m_n))$, par la
projection canonique, n'est pas bornŽe.  La projection de
$(H_{M,\sigma}(m_n))$ sur $\a_{G,\sigma}$ paralllement ˆ $\a^{G}\cap \a_{M, \sigma}$ est
Žgale ˆ
$(H_{G,\sigma}(g_n))$, car $H_{G,\sigma}(h_n)=0$ et $H_{G,\sigma}(u_n)= H_{G}(u_n)=0$. C'est
une suite dans un rŽseau, qui est convergente d'aprs l'hypothse sur $(g_n)$. Quitte ˆ
extraire une sous-suite on peut supposer que cette suite est constante, ce que l'on fait. Il
faut dŽmontrer que
$(H_{M,\sigma}(m_n))$ n'est pas bornŽe.\ste  Supposons
$(H_{M,\sigma}(m_n))$ bornŽe. Comme $H_{M,\sigma}$ prend ses valeurs dans un
sous-groupe discret de $\a_{M,\sigma}$, quitte ˆ extraire une sous-suite on peut supposer
que $H_{M,\sigma}(m_n)$ est constant. Quitte ˆ multiplier, ˆ droite, 
$(g_n)$ par un ŽlŽment de $M$ convenable, on dispose d'une suite $(g_n)$ convergeant vers
$g \in \overline{HP}\backslash HP$ avec $g_n = h_nm_nu_n$ et $H_{M,\sigma}(m_n) =
0$, i.e. $m_n \in M^{1}_{\sigma}$.\ste 
En utilisant la Proposition   \ref{pitilde} (ii), on voit que: $${\tilde
\pi}_i(m_n){\tilde v}_i = {\tilde v}_i,
\>\> n
\in
\N, i = 1,...,m'$$ Par ailleurs $({\sigma(g_n)^{-1}}g_n)$ converge. Donc, posant
$\overline{u_n} =
\sigma(u_n)^{-1} \in \overline{U}$, o $\overline{U} = \sigma(U)$, on a
$(\overline{u_n}(\sigma(m_n^{-1})m_nu_n)$ qui converge.\ste 
On montre un analogue du Lemme 5 de [BrD] qui permet de voir que $\overline{u_n}$
converge {de mme que $u_n$. Il en va donc de mme de $\sigma(m_n)^{-1}m_n$.\ste  
Notons $G^{\sigma}$  (resp. $M^{\sigma}$) le groupe des points fixes de $\sigma$ dans $G$
(resp. $M$).\ste 
L'application de $M^\sigma \backslash M$ dans $M$,  $M^\sigma m \mapsto
\sigma(m)^{-1} m$, est un homŽomorphisme de
$M^\sigma \backslash M$ sur son image. Ceci rŽsulte du Lemme \ref{orbites} de l' Appendice,
appliquŽ ˆ l'action par conjugaison tordue de $M$ sur $M^{-\sigma} = \{ m \in M \vert 
\sigma(m) = m^{-1} \}$:$$m.x = mx\sigma(m)^{-1},m,x \in M$$
En effet $M$ n'a qu'un nombre fini d'orbites dans $M^{-\sigma}$ (cf. [R]).
 \ste 
On en dŽduit que $({M^\sigma}m_n)$ converge, ce qui veut dire qu'il existe une suite
$(m'_n)$ dans $M^\sigma$ telle que $(m'_nm_n)$ converge. Mais $
M^\sigma/ ({M \cap H}) $ est fini. Donc, quitte ˆ extraire une sous-suite, on peut trouver
$(h'_n)$, suite dans $M \cap H$, telle que $(h'_nm_n)$ converge dans $M$. On Žcrit alors $g_n =
h_n {h'_n}^{-1}h'_n m_n u_n$  i.e.  ${g_n} = {h''_n} {m''_n} {u_n}$ avec ${h''_n} = {h_n}
{h'_n}^{-1},\>\> {m''_n} = {h'_n} {m_n}$.\ste 
Alors $(m''_n)$ converge dans $M$ de mme que $(u_n)$ dans $U$ donc aussi $(h''_n)$ dans
$H$. Mais alors on aurait $g \in HP$. Une contradiction qui montre que l'image de
$(H_{M,\sigma}(m_n))$ dans $\a_{M,\sigma} / \a_{G,\sigma}$ est non bornŽe.
Mais pour $i = 1,...,m'$: $$
\varepsilon_i(g_n) = e^{n_i\tilde{\delta}_i(H_{M,\sigma}(m_n))}$$  et
$(\varepsilon_i(g_n))$ converge vers $\varepsilon_i(g)$. Si pour tout $i = 1,...,m'$,
$(\varepsilon_i(g_n))$ convergeait vers une limite non nulle, on en dŽduirait que
$(H_{M,\sigma}(m_n))$ serait bornŽe modulo $\a_{G,\sigma}$, car
$({\tilde \delta}_1,...,{\tilde \delta}_{m'})$ est une base de ${(\a_{M,\sigma} /
\a_{G,\sigma})}^*$. Donc l'un des $\varepsilon_i(g)$ est nul, comme dŽsirŽ.
\ste (iii) rŽsulte de la preuve de (ii).\qed\ste
Si $ i \in \{1,...,l\}$, on fixe une base de $V_i$ formŽe de vecteurs poids sous
$A_0$, ce qui permet de dŽfinir une norme sur  $V_i$ en prenant le maximum des valuations
des coordonnŽes dans cette base. Puis on en  dŽduit une norme sur $\tilde{V}_i^*=End V_{i}$,
en prenant le maximum des valuations des coefficients de la matrice dans cette base.
\ste  On pose, pour $g \in G$: 
\ber $\Vert gH\Vert_i = \Vert\tilde{\pi}_i^*(g) {\tilde e}_{i,H}^*\Vert$, \space $i =
1,...,l$
\ste 
$\Vert gH\Vert_0 = e^{\Vert H_{G,\sigma}(g)\Vert}$ \eer
o l'on a muni $\a_{G,\sigma}$ de la norme provenant du produit scalaire sur $\a_0$.\ste 
On dŽfinit: \beq \Vert gH\Vert = \prod_{i=o}^l \Vert gH\Vert_i,\>\> g \in G   \eeq
\begin{rem}On remarque que $\tilde{\pi}_i^*(g) {\tilde e}_{i,H}^*$ est Žgal ˆ
$\pi_i(g^{\sigma}g^{-1})$ et $\Vert gH\Vert_i$ est la norme de l'opŽrateur  $\pi_i
(g^{\sigma}g^{-1})$ associŽe ˆ la norme introduite ci-dessus sur $V_i$.

\end{rem}
On dŽfinit de manire similaire une application $m \to \Vert m (M\cap H) \Vert$ pour 
$m \in M$, en utilisant les reprŽsentations rationnelles de $M$ associŽes aux
$n_i\tilde{\delta}_i^M$, $i \in \{m+1,...,l\}$ (cf. Proposition \ref{deltaM}).  
\begin{prop}: Il existe $\nu_0 \in \a_{M,\sigma}^*$,  $P$-dominant, tel que
pour toute suite
$(g_n)$ dans $HP$, convergeant vers un ŽlŽment \space $g$ de $\overline{HP}$
$\backslash$ $HP$, on ait la propriŽtŽ suivante:\ste Ecrivant \space  $g_n = h_nm_nu_n$
avec $h_n \in H, m_n \in M , u_n \in U$, il existe $C >0$ tel que la suite $\Vert
m_n^{-1}(M \cap H)\Vert$ soit bornŽe par
$C e^{-\nu_0(H_{M,\sigma}(m_n))}$.\ste 
 \end{prop} 
\dem  On choisit $i \in \{m+1,...,l\}$ i.e. $\alpha_i \in \Theta$ et une base de
$\tilde{V}_i^M$, $(\tilde{v}_{i,j})$, construite ˆ partir d'une base de  vecteurs poids sous
$A_0$ de $V_i^M$. On pose:$$\varepsilon_{i,j}(g): = \vert <\tilde{\pi}_i(g) \tilde{v}_{i,j},
\tilde{e}_{i,H}^*   >\vert_{\F},
\>  \> g\in G$$ Alors $((\varepsilon_{i,j})(g_n))$ est bornŽe car elle converge vers
$((\varepsilon_{i,j})(g))$.
 De plus,  gr\^ace  ˆ la Proposition \ref{deltaM} (ii),  on voit que: 
$$\varepsilon_{i,j}(g_n) = \vert{\tilde \Lambda}_{i,M}(m_n)\vert_{\F}\>\>\vert
<(\tilde{\pi}_i^M)(m_n) {\tilde v}_{ij}, 
\tilde{e}_{i,M \cap H}^* >\vert_{\F} =
e^{n_i\tilde{\delta}_i(H_{M,\sigma}(m_n))} \vert <(\tilde{\pi}_i^M)(m_n) {\tilde v}_{ij}, 
\tilde{e}_{i,M \cap H}^* >\vert_{\F}$$
En passant au sup sur $j$ et en posant $\tilde{\nu}_i = n_i
{{\tilde \delta}_i}{}_{\vert{\a_{M,\sigma}}}$, on trouve: $$Sup_j  \varepsilon_{i,j}
(g_n) = e^{\tilde{\nu}_i({H_{M,\sigma}} (m_n))} \Vert m_n^{-1} (M \cap H) \Vert _{i}, i=m+1,
\dots,l $$ On dŽduit de ce qui prŽcde que: \ber \label{normei} $\prod_{i=m+1}^l \Vert m_n^{-1}
(M
\cap H )\Vert_i$ est bornŽe par un multiple de $\prod_{i=m+1}^l 
e^{-\tilde{\nu}_i (
H_{M,\sigma} (m_n))}$ pour $ i=m+1,
\dots,l$ \eer Par ailleurs, pour $i=1,...,m'$, $(\varepsilon_i(g_n))$ est
bornŽe.\ste  Mais  $\varepsilon_i(g_n)$ =
$e^{n_i {\tilde \delta}_i(H_{M,\sigma}(m_n))}$, donc, pour $i=1,...,m'$,
$(n_i\tilde{\delta}_i({H}_{M,\sigma} (m_n)))$ est bornŽe supŽrieurement  de sorte
que:\ber \label{norme0} Pour $i=1,...,m'$, $(e^{\vert
n_i\tilde{\delta}_i({H}_{M,\sigma} (m_n) \vert})$ est bornŽe par un multiple de
$(e^{-{n_i} \tilde{\delta}_i({H}_{M,\sigma} (m_n)})$. 
\eer
Revenant ˆ la dŽfinition de $\Vert m(M\cap H)\Vert_0 $, on voit qu'il existe $c>0$ tel que: $$
\Vert m(M\cap H)\Vert_0\leq  e^{(c \sum_{i=1, \dots, m'}\vert
n_i\tilde{\delta}_i({H}_{M,\sigma} (m))\vert )+\Vert  H_{G, \sigma}(m)\Vert}, m \in M $$
On notera $\tilde{\nu}_0 = \sum_{i=1}^{m'}c n_i \delta_i$ et $\nu_0 = \tilde{\nu}_0 +
\sum_{i=m+1}^l \tilde{\nu}_i$. Alors d'aprs (\ref{norme0}) et (\ref{normei}) et le fait que
$(H_{G, \sigma}(m_n))= (H_{G, \sigma} (g_n))$ est bornŽe,   on voit que
$\nu_0$, vŽrifie les propriŽtŽs voulues.\qed
\subsection{Une autre prŽsentation des fonctions j}Les racines de $A_M$ dans l'algbre
de Lie de $G$ s'identifie ˆ des ŽlŽments de $\a_M^{*}$.  Soit   $\rho_P$
 la demi-somme des racines de $A_M$ dans l'algbre de Lie de $U$ (ou $P$), comptŽes
avec les multiplicitŽs.\ste On note
$C(G,P,-2\rho_P)$ l'espace des fonctions continues, $f:G\to\C$,   telles que:$$f(gmu) =
e^{-2\rho_P({H}_M (m))} f(g)  ,\>\> g \in G , m \in M , u \in U $$   Alors: \ber
\label{forme invariante} Si $f\in C(G,P,-2\rho_P)$, la restriction de $f$ ˆ
$K$ est invariante ˆ droite par $K\cap P$ et la forme linŽaire sur $C(G,P,-2\rho_P)$:
$f
\to
\int_{K/K \cap P} f(k(K
\cap P)) d\dot{k} $ est invariante par  translation ˆ gauche par les  ŽlŽments de $G$.
\eer 
Soit $C(G,P,\delta^*,\chi)$ l'espace des fonctions $\psi: G \to V_\delta^*$, o
$V_\delta^*$ est le dual algŽbrique de $V_{\delta}$,  qui sont faiblement continues i.e.
telles que pour tout $v \in V_\delta$, $g \mapsto <\psi(g) , v >$ est continue sur $G$
et telles que \beq \psi(gmu) = e^{-2\rho_P({H}_M (m))} \chi(m) {\delta^*}(m)^{-1}
\psi(g) \eeq  
\begin{rem} Le facteur $2$ dans $2\rho_P$ tient compte du fait que nous avons choisi
l'induction non  normalisŽe. \end{rem} \ste Alors si $\psi \in C(G,P,\delta^*,\chi)$ et
$\varphi \in I_{\chi}= ind_P^G V_{\delta \otimes \chi}$, l'application $g \mapsto
<\psi(g),\varphi(g)>$ est continue sur $G$, car $\varphi$ est localement constante,  et
c'est un ŽlŽment de $C(G,P,-2\rho_P)$.\ste 
On pose:$$<\psi,\varphi> = \int_{K/K \cap P} <\psi(k) , \varphi(k)>
d\dot{k}$$ Ceci dŽfinit un crochet de dualitŽ $G$-invariant entre $C(G,P,\delta^*,\chi)$
et $I_{\chi}$ qui permet de regarder $C(G,P,\delta^*,\chi)$
comme un sous-espace de $I_{\chi}^*.$ 

\begin{theo} Soit $(\delta,V_\delta)$ une reprŽsentation admissible de type fini $M$,  $\eta
\in V_\delta^{*M \cap H}$ et  $r \in \R$ tels que: \ste Pour tout $v \in V_\delta$, il
existe
$C>0$  tels que: $$\vert <m \eta , v >\vert \leq C \Vert m(M\cap H) \Vert^r ,
\>\> m\in M$$ Soit $\nu_{0}$ comme dans la Proposition 13.  
Soit $\chi \in X(M)_{\sigma}$  et vŽrifiant $Re\chi -2\rho_P-
r\nu_0$ strictement $P$-dominant. \ste 
On dŽfinit une application  $\varepsilon(G,P,\delta,\chi, \eta)$ de $  G$ dans  $V_\delta^*$
par:
$$\varepsilon(G,P,\delta,\chi, \eta)(hmu) = e^{-2\rho_P(H_M (m))} \chi(m) 
\delta^*(m^{-1}) \eta , h \in H , m \in M , u \in U$$ $$\varepsilon(G,P,\delta,\chi, \eta
)(g) =0
\>\>   si \>\> g \notin HP$$
 \ste (i) Alors  $\varepsilon(G, P,\delta,\chi, \eta)$ est un ŽlŽment de
$C(G,P,\delta^*,\chi)$, 
$H$-invariant ˆ gauche.
\ste(ii) Soit $q\in B$ comme dans le ThŽorme 1 (i). On suppose  que $\chi$ vŽrifie les
hypothses de (i) et que $q(\chi)$ est non nul. \ste L'ŽlŽment de
$I_{\chi}^{*H}$ correspondant ˆ l'ŽlŽment $\varepsilon(G,P,\delta,\chi, \eta
) $ de $C(G,P,\delta^*,\chi)$, est Žgal ˆ $j(P,\delta,\chi,\eta)$.
\end{theo}
\dem La partie (i) est une consŽquence immŽdiate des propriŽtŽs de $\eta$,  de la
Proposition 12 (iii) et de la Proposition 13.\ste 
Prouvons (ii). On a l'analogue de la formule du Lemme 1.3 de [Ol] pour les groupes
$p$-adiques (c.f. l.c., ThŽorme 7.1 ) qui implique:
\ber  Soit $f$ une fonction continue sur $K/K \cap P$ nulle en dehors de $(HP)\cap K$,
alors, pour une normalisation convenable des mesures:$$\int_{K/K \cap P} f(k(K \cap P))
d\dot{k} = \int_{H/H \cap M} f(k(h)(K \cap P)) e^{-2\rho_P(H_M(h))} d\dot{h}$$o l'on a
Žcrit $h = k(h)p$ avec $k(h) \in K$ et $p \in P$. La classe de $k(h)$ modulo $K \cap P$ est
bien dŽfinie. De plus $\rho_P$ Žtant antiinvariante par $\sigma$, $h \mapsto
\rho_P(H_M(h))$ est invariante ˆ droite par $H \cap P = H \cap M$. \eer
Si $ f \in C(G,P,-2\rho_P)$, sa restriction ˆ $K$ est $K \cap P$ invariante et on
a:$$f(h) = f(k(h)(K \cap P)) e^{-2\rho_P(H_M(h))} , \>\> h \in H$$ De plus la
restriction de
$f$ ˆ
$H$ est
$H
\cap M$-invariante.\ste
Alors il rŽsulte de ce qui prŽcde que:
\ber Pour $f \in C(G,P,-2\rho_P)$ nulle en dehors de $HP$, on a: $$\int_{K/K \cap P}
f(k(K
\cap P)) d\dot{k} = \int_{H/H \cap M} f(h) d\dot{h}$$ \eer
Soit $\varphi \in I_\chi$ avec $\chi$ comme dans (ii). La dŽfinition de
$\varepsilon(G,P,\delta,\chi,\eta)$ et son invariance ˆ droite par $H$, jointe ˆ l'Žquation
prŽcŽdente montre que: \beq \int_{K/K \cap P}<\varepsilon(G,P,\delta,\chi,\eta)(k(K \cap P)),
\varphi(k(K \cap P)) >d\dot{k}= \int_{H/H \cap M}<\eta, \varphi(h) > d\dot{ h}  \eeq 
Ceci implique en particulier que la restriction de l'ŽlŽment $\varepsilon_\chi$ de
$I_\chi^{*H}$ ˆ $J_\chi$, correspondant ˆ $\varepsilon(G,P,\delta,\chi,\eta)$ est Žgal, avec
les notations du ThŽorme 2 (ii),  ˆ $\xi(P,\delta,\chi,\eta) \in J_\chi^{*H}$.\ste
Alors, d'aprs le ThŽorme 2 (iii), cet ŽlŽment est Žgal ˆ $j(P,\delta,\chi,\eta)$. \qed
\begin{rem} Si $P$ est Žgal ˆ $P_\emptyset$ la condition sur $\eta$ est toujours
satisfaite car d'aprs [HWan], Lemme 4.5 ( cf. aussi [HH],  et 
Lemme 1.9),
$M_\emptyset/M_\emptyset
\cap H$ est compact modulo le centre de $M_\emptyset$. De plus, si  $A_{0}$ est un tore
dŽployŽ maximal de $M_{\emptyset}$, il est $\sigma$-stable et $A_{\emptyset}$ est
l'unique tore
$\sigma
$-dŽployŽ maximal  de $M_{\emptyset}$ (cf. l.c.). Tenant compte de
[HH],  1.3, on en dŽduit que $M_\emptyset/M_\emptyset
\cap H$ est compact modulo $A_{\emptyset}$. \ste Des rŽsultats rŽcents de Nathalie Lagier,
indiquent  que, pour n'importe quel $\sigma$-sous-groupe parabolique,  la condition sur
$\eta$ est toujours   satisfaite,  au moins si $\delta$ est unitaire.  

\end{rem}
 \section{ Appendice:   Rappels sur les l-groupes }
\begin{lem}\label{orbites}
  (i) Soit $G$ un $l$-groupe, $X$ un $l$-espace, i.e. un 
espace
topologique avec une base d'ouverts formŽe d'ensembles compacts.
  On suppose que $X$ est muni d'une action continue de $G$, i.e.
que l'application
 $G\times X
 \rightarrow X$ est continue,  que $G$ n'a qu'un nombre fini d'orbites dans
$X$  et que $G$ est
 dŽnombrable ˆ l'infini.  Alors: \ste  (i) Il existe une orbite ouverte.\ste
                                 (ii) Il existe des ouverts $G$-invariants,
 $U_{0} =
 \emptyset\subset U_{1}\dots \subset U_{n} = X$, o pour $i\geq 0$,
 $U_{i+1} \setminus U_{i}$
 est une
 $G$-orbite.  En particulier les orbites sont localement fermŽes, donc sont des
$l$-espaces, i.e.  admettent une base de voisinages ouverts et compacts. \ste 
                                 (iii) Si $X_{1}$ est une orbite $G{x_1}$ et
$H$ est le  stabilisateur de $x_{1}$, alors $H$ est un sous groupe fermŽ de
$G$ et l'application: $G/H 
\rightarrow X_{1}$ donnŽe par $gH \rightarrow gx_{1}$ est un homŽomorphisme. 
\end{lem}
\dem 
Prouvons (i). Soit $x_1, \dots , x_p $ des  reprŽsentants des orbites sous $G$
dans $X$, qui est fini par hypothse. Soit ($K_{n}$) un ensemble dŽnombrable d'ouverts compacts
recouvrant
$G$.  Si $X = \bigcup_{i,n}K_{n}x_{i}$, les $K_{n}x_{i}$ sont compacts donc
fermŽs. D'aprs le thŽorme de Baire, l'un de ces ensembles, $K_{n}x_{i}$, est
d'intŽrieur non vide, donc
$G{x_i}$
 est ouvert 
dans $X$.\ste 
Prouvons (ii).  On construit $U_{i}$ par rŽcurrence sur $i$. On prend pour
$U_{i}$ une orbite  ouverte et pour $U_{i+1}$ la rŽunion de $U_{i}$ et d'une
orbite ouverte dans $X
\setminus U_{i}$.\ste  
Prouvons (iii).  L'application $p: G/H \rightarrow X_{1}$ est
clairement continue et bijective.  Il reste ˆ voir que cette application est
ouverte. En utilisant les translations, il  suffit pour cela de voir que tout
voisinage de $eH$ dans $G/H$ s'envoie sur un  voisinage de $x_1$ dans $X_1$.
Tout voisinage ouvert compact $V$ de $eH$ dans
$G/H$   contient l'image  par la projection naturelle de $G$ dans $G/H$, d'un sous groupe
ouvert compact $K$. Montrons que $G = \bigcup_{n \in N}g_{n}K$ pour une suite $(g_{n})$
d'ŽlŽments de $G$.  En effet $G = \bigcup_{n}K_{n}$ et chaque $K_{n}$ Žtant compact,
il est recouvert  par un nombre fini d'ensembles de la forme $gK$. Alors les
$(g_{n}K{x_1})_{n \in \N}$ sont fermŽs et recouvrent $X_{1}$. Donc, toujours
par le thŽorme de Baire, l'un d'eux, $g_{n}Kx_{1}$, est d'intŽrieur non vide.
La translation
 $x \rightarrow   g_{n}^{-1}x$ Žtant un  homomŽorphisme, on en dŽduit que $Kx_{1}$ est
d'intŽrieur non vide. Quitte ˆ utiliser encore une translation, par un ŽlŽment de $K$ cette
fois, on voit que $K{x_1}$  contient un voisinage de $x_1$. Ainsi $p(V)$
contient un voisinage de $x_1$ comme dŽsirŽ. Ceci achve de prouver (iii).\qed
 Rappelons un rŽsultat de [M] ( Corollaire 2): \ber Soit $G$ un $l$-groupe,
$H$ un sous-groupe fermŽ, le fibrŽ principal: $G \rightarrow G/H$ est trivial, i.e. il
existe une section continue $s: G/H \rightarrow G$ de $p$.\eer

\begin{lem} \label{fibre}
L'application linŽaire $T:  C_c^\infty(G/H) \otimes C_c^\infty(H)
\rightarrow C_c^\infty(G)$, dŽfinie par: $(T(\psi \otimes \eta
))(s(x)h) = \psi(x)\eta(h)$, est  un entrelacement bijectif entre le
produit tensoriel de la reprŽsentation triviale de $H$ sur
$C_c^\infty(G/H)$ avec la reprŽsentation rŽgulire droite sur
$C_c^\infty(H)$ et la reprŽsentation rŽgulire droite de $H$ sur
$C_c^\infty(G)$. De plus l'inverse de $T$, $T^{-1}$ est donnŽ par:
$(T^{-1}\varphi)(x,h) =
\varphi(s(x)h)$.
  On a un rŽsultat analogue pour l'action rŽgulire gauche.\end{lem}
\dem  D'abord $T$ est bien ˆ valeurs dans $C_c^\infty(G)$, car, pour un
$l$-espace $X$,  un ŽlŽment de
$C_c^\infty(X)$ est exactement une application continue sur $X$ prenant un nombre
fini de valeurs et ˆ support compact. De mme on voit que si $\varphi \in
C_c^\infty(G)$, l'application: $((x,h) \mapsto \varphi(s(x)h)$, est
ŽlŽment de
$C_c^\infty(G/H \times H) \cong C_c^\infty(G/H) \otimes C_c^\infty(H)$, que l'on
note $T'\varphi$. Clairement $T$ et $T'$ sont des applicatons inverses l'une
de l'autre. D'autre part $T$ ˆ la propriŽtŽ d'entrelacement voulue. D'o le
Lemme.   \qed

\section{Bibliographie}
\parindent=0pt
\hfuzz=30pt

\noindent[BeZ] Bernstein, I. N., Zelevinsky, A. V., Induced representations of reductive 
$\p$-adic groups. I , Annales Scientifiques de l'ƒcole Normale SupŽrieure SŽr. 4, 10 (1977),
Ê441-472. \ste

\noindent[Bl] Blanc, P., Projectifs dans la catŽgorie des G-modules topologiques.
C.R.Acad.Sc.Paris, t.289 (1979), 161-163.\ste

\noindent[BlBr] Blanc, P., Brylinski, J-L,.,  Cyclic homology and the
Selberg principle.  J. Funct. Anal.  109  (1992) 289--330.\ste

\noindent [BoT] Borel, A, Tits, J., Groupes rŽductifs, Publications MathŽmatiques de
l'IHES, 27(1965) 55--151.\ste

\noindent[BoWall] Borel, A., Wallach N., Continuous cohomology, discrete
subgroups, and representations of reductive groups. Second edition.
Mathematical Surveys and Monographs,
 67. American Mathematical Society, Providence, RI, 2000. \ste

\noindent[BrD] Brylinski, J.L., Delorme,  P., Vecteurs distributions $H$-invariants pour les
sŽries principales gŽnŽralisŽes d'espaces symŽtriques rŽductifs et prolongement mŽromorphe
d'intŽgrales d'Eisenstein.   Invent. Math.  109  (1992)619--664.\ste

\noindent[CarD] Carmona, J., Delorme, P., Base mŽromorphe de vecteurs
distributions $H$-invariants pour les sŽries principales gŽnŽralisŽes d'espaces symŽtriques
rŽductifs: Žquation fonctionnelle. J. Funct. Anal.  122  (1994) 152--221.\ste

\noindent [CartE], Cartan, H., Eilenberg,  S., Homological algebra, Princeton Univ. Press,
1956. \ste 

\noindent [Cas] Casselman, W. A new nonunitarity argument for $p$-adic
representations.  J. Fac. Sci. Univ. Tokyo,  28  (1982) 907--928.\ste

\noindent[D], Delorme, P.,  Harmonic analysis on real reductive symmetric spaces. Proceedings of the
International Congress of Mathematicians, Vol. II  (Beijing, 2002), 545--554, Higher
Ed. Press, Beijing, 2002. \ste

\noindent[G] Guichardet, A., Cohomologie des groupes topologiques et des
algbres de Lie.  CEDIC, Paris, 1980. \ste

\noindent[HH]  Helminck, A. G.,  Helminck, G. F.,  A class of parabolic $k$-subgroups
associated with symmetric $k$-varieties.  Trans. Amer. Math. Soc.  350  (1998)
4669--4691.\ste

\noindent[HWan] Helminck, A. G.,  Wang, S. P. On rationality properties of involutions of
reductive groups.  Adv. Math.  99  (1993) 26--96.\ste

\noindent[Hi] Y. Hironaka, Spherical functions and local densities on Hermitian forms, J.
Math Soc. Japan, 51 (1999) 553-581.  \ste 

\noindent[HiSat] Y. Hironaka, F. Sato, Spherical functions and local densities of alternating
forms, Am. J. Math., 110 (1988), 473-512. \ste 

\noindent[Hu1] Humphreys, J.E. , Introduction to Lie algebras and representation theory. Second
printing, revised. Graduate Texts in Mathematics, 9. Springer-Verlag, New York-Berlin, 1978.
\ste 

\noindent[Hu2] Humphreys, J.E.,  Linear algebraic groups. Graduate Texts in Mathematics, No.
21. Springer-Verlag, New York-Heidelberg, 1975. \ste

\noindent[M] Michael E. Selected selection theorems, Amer. Math. Monthly, 63 (1956),
223-238.\ste 

\noindent[O] Offen, O.,  Relative spherical functions on $\wp$-adic symmetric spaces (three 
cases). Pacific J. Math. 215 (2004), no. 1, 97--149.\ste 

\noindent[OS] Offen, O., Sayag, E., On unitary distinguished representations of $GL_{2n}$
distinguished by the symplectic group. preprint \ste 

\noindent[Ol] Olafsson, G., Fourier and Poisson transformation associated to a semisimple
symmetric space.  Invent. Math.  90  (1987), 605--629. \ste

\noindent[R] Richardson, R.W., Orbits, invariants and representations associated to
involutions of reductive groups. Invent Math. 66 (1982), 287-312.\ste
  
\noindent[Wald] Waldspurger, J.-L. La formule de Plancherel pour les groupes $p$-adiques
(d'aprs Harish-Chandra). J. Inst. Math. Jussieu  2  (2003) 235--333.

}\ste

\ste \ste 
Philippe Blanc \ste 
e-mail: blanc@iml.univ-mrs.fr\ste 
Patrick Delorme \ste e-mail: delorme@iml.univ-mrs.fr\ste \ste 

\ste Institut de MathŽmatiques de Luminy, UMR 6206 CNRS
\ste UniversitŽ de la MŽditerranŽe
163 Avenue de Luminy\ste Case 907\ste 
13288 Marseille Cedex 09
France
\end{document}